\documentclass{elsart}
\usepackage{natbib,amssymb}

\newenvironment{remark}{\refstepcounter{thm}
 \bigbreak\noindent{\bf Remark \arabic{thm} }}
 {\medbreak}

\newenvironment{definition}{\refstepcounter{thm}
 \bigbreak\noindent{\bf Definition \arabic{thm} }}
 {\medbreak}

\newenvironment{example}{\refstepcounter{thm}
\bigbreak\noindent{\bf Example \arabic{thm} }}
 {\medbreak}

\newcounter{itemcounter}
\newenvironment{items}{
   \begin{list}{\alph{itemcounter})}
   {\usecounter{itemcounter}\setlength{\topsep}{3 pt}
   \setlength{\partopsep}{0 pt}\setlength{\itemsep}{0 pt}
      \setlength{\labelwidth}{2 em}
   }}{\end{list}}

\newenvironment{rules}{
   \begin{list}{\alph{itemcounter})}
   {\usecounter{itemcounter}\setlength{\topsep}{3 pt}
   \setlength{\partopsep}{0 pt}\setlength{\itemsep}{0 pt}
      \setlength{\labelwidth}{1.65 em}
   }}{\end{list}}

\newenvironment{rulesp}{
   \begin{list}{\alph{itemcounter})}
   {\usecounter{itemcounter}\setlength{\topsep}{3 pt}
   \setlength{\partopsep}{0 pt}\setlength{\itemsep}{0 pt}
      \setlength{\labelwidth}{1.7 em}
   }}{\end{list}}

\def\cocoa
 {\mbox{\rm C\kern-.13em o\kern-.07em C\kern-.13em o\kern-.15em A}}

   \font\tengothic=eufm10
   \font\sevengothic=eufm7
   \font\fivegothic=eufm5

   \font\tenmsy=msbm10
   \font\sevenmsy=msbm7
   \font\fivemsy=msbm5

\newfam\msyfam
      \textfont\msyfam=\tenmsy
      \scriptfont\msyfam=\sevenmsy
      \scriptscriptfont\msyfam=\fivemsy
\newfam\gothicfam
   \textfont\gothicfam=\tengothic
   \scriptfont\gothicfam=\sevengothic
   \scriptscriptfont\gothicfam=\fivegothic

\def\Cal#1{{\mathcal #1}}

\let\sem=\bf
\let\phi=\varphi
\let\rho=\varrho
\let\theta=\vartheta
\let\epsilon=\varepsilon

\def\TTo#1{\mathop{\longrightarrow}\limits ^{#1}}

\def\squareforqed{\hbox{\rlap
{$\sqcap$}$\sqcup$}}
\def\qed{\ifmmode\squareforqed\else
{\unskip\nobreak\hfil
\penalty50\hskip1em\null\nobreak\hfil
\squareforqed
\parfillskip=0pt\finalhyphendemerits=0\endgraf}
\fi}

\def\proof{\rm \trivlist
 \item[\hskip \labelsep{\it Proof.}]}
\def\endproof{\qed \endtrivlist}

\def\NF{\mathop{\rm NF}\nolimits}
\def\NR{\mathop{\rm NR}\nolimits}
\def\LT{\mathop{\rm LT}\nolimits}

\def\LM{\mathop{\rm LM}\nolimits}

\def\Mat{\mathop{\rm Mat}\nolimits}

\def\rk{\mathop{\rm rk}\nolimits}

\def\dsum{\mathop{\oplus}\nolimits}

\def\iff{\Leftrightarrow}
\def\Tnr{\mathbb{T}^n\langle e_1,\dots,e_r\rangle}
\def\tfrac#1#2{\textstyle{#1\over #2}}
\def\lcm{\mathop{\rm lcm}\nolimits}
\def\Syz{\mathop{\rm Syz}\nolimits}

\journal{Journal of Symbolic Computation}
\begin{document}

\begin{frontmatter}
\title{Efficiently Computing\\ Minimal Sets of Critical Pairs}

\author[Genova]{M. Caboara\corauthref{massimo}},
\corauth[massimo]{Corresponding author.}
\ead{caboara@dima.unige.it}
\author[Dortmund]{M. Kreuzer},
\ead{Martin.Kreuzer@mathematik.uni-dortmund.de}
\author[Genova]{L. Robbiano}
\ead{robbiano@dima.unige.it}

\address[Genova]{Department of Mathematics, University of Genoa, Italy}
\address[Dortmund]{Fachbereich Mathematik, Universit\"at Dortmund, Germany}

\begin{abstract}
In the computation of a Gr\"obner basis using Buchberger's
algorithm, a key issue for improving the efficiency is to
produce techniques for avoiding as many unnecessary
critical pairs as possible. A good solution would be to
avoid {\it all}\/ non-minimal critical pairs, and hence to
process only a {\it minimal set of generators}\/ of the module
generated by the critical syzygies. In this paper we show
how to obtain  that desired solution in the homogeneous case
while retaining the same efficiency as with the classical implementation.
As a consequence, we get a new Optimized Buchberger Algorithm.
\end{abstract}

\end{frontmatter}


\section{Introduction}
\label{Introduction}

Ever since practical implementations of Buchberger's famous
algorithm for computing Gr\"obner bases became
feasible \citep{B1965}, it has been clear that, in order to improve the
efficiency of this algorithm, one needs to avoid the
treatment of as many critical pairs as possible.
\citet{B1979} studied this problem for the
first time, and later in~\citep{B1985} and~\citep{GM1987}
his results were substantially improved and expanded.
Never\-the\-less, \citet{GM1987} showed that
their method did not always produce a minimal set of generators
of the module generated by the critical syzygies.
However, their method was very efficient and yielded an
{\it almost} minimal set of critical pairs. Since then,
many kinds of optimizations of Buchberger's algorithm have
been found, in particular by implementers of computer algebra
systems. But the problem of efficiently minimalizing the
critical pairs has gone largely unnoticed and seems to be overdue
for a solution. Indeed, that is the main objective of this paper.

To achieve our goal, we proceed as follows. Foremost, we
need a detailed understanding of the entire process of
computing Gr\"obner bases, in particular in the homogeneous
case. An algorithm for simultaneously computing a
Gr\"obner basis and a minimal system of generators contained
in it is fine-tuned when the input is a reduced Gr\"obner basis.
Then this result is applied to critical
syzygies, using the fact that we show how the old criteria
$M(i,j)$ and $F(i,j)$ of~\citep{GM1987}
yield a reduced Gr\"obner basis of the module of syzygies of the 
leading terms. Besides, when applied to this special case, the
algorithm admits many subtle optimizations. In the end, we
really achieve the goal of minimalizing the critical pairs
efficiently.

\medskip
\centerline{\it Now, why do we think that what we achieved
is important?}
\smallskip

The first reason is theoretical curiosity. It is common knowledge
among the implementers of Buchberger's algorithm that the
criteria of  Gebauer and M\"oller {\it almost}\/ produce a minimal
set of critical pairs. We wanted to see whether that {\it vox
populi}\/ is really true. Of course one could use a standard
minimalization process to produce minimal sets of critical pairs,
but this method could only handle small examples. Instead, we
observed that, after applying two of the criteria of Gebauer and
M\"oller, a {\it reduced Gr\"obner basis}\/ of the module of
syzygies of the leading terms is obtained. 
Then we were able to see the difference
between the reduced Gr\"obner basis and a minimal set of generators
of this module, and how this difference depends on the size
of the example.

Another important reason is that we wanted to be able to
compute a minimal set of generators of this module
with the {\it same efficiency}\/ as the usual application of the
Gebauer-M\"oller criteria. And we wanted to do it
while computing a Gr\"obner basis, so that we can
replace the Gebauer-M\"oller criteria by our procedure.
As we show in the last sections, we achieved this goal.

A third reason is that our results hold in full generality,
namely for Gr\"obner bases of modules over positively
(multi-) graded rings. Other optimizations of
Buchberger's algorithm, e.g.\ ideas using
trivial syzygies (see for instance \citet{F2002}), 
do not hold in this generality.
Moreover, we would like to point out that the pairs
we discard are truly useless, whereas pairs
between elements in a reduced Gr\"obner bases which reduce
to zero can still be useful for the computation of syzygies.

Finally, the readers should know that the basic terminology
is taken from the book of the second and third
authors\,\citep{KR2000}.

\section{Some Background Material}
\label{Some Background Material}

Since we are interested in optimizing Buchberger's algorithm in
the homogeneous case, we start by saying which gradings we
consider. From now on let~$K$ be a field and
$P=K[x_1,\dots,x_n]$ a polynomial ring over~$K$. Moreover, let
$W\in\Mat_{m,n}(\mathbb{Z})$ be an $m\times n$-matrix with integer
entries. Then there exists exactly one $\mathbb{Z}^m$-grading on~$P$
such that every term $t=x_1^{\alpha_1}\cdots x_n^{\alpha_n}$ is
homogeneous of degree $\deg_W(t)=W\cdot
(\alpha_1,\dots,\alpha_n)^{\rm tr}$. We say that~$P$ is {\sem
(multi\hbox{-)} graded} by~$W$. The matrix~$W$ is called the {\sem
degree matrix} and its rows are called the {\sem weight vectors}.

For instance, the grading on~$P$ given by~$W=(1,\dots,1)$
is the standard grading. For every $d\in \mathbb{Z}^m\!$,
the homogeneous component of degree~$d$ of~$P$ is
$P_{W,d}= \dsum_{\deg_W(t)=d}K\cdot t$.
Given $\delta_1,\dots,\delta_r\in \mathbb{Z}^m$, the graded
free $P$-module $F=\dsum_{i=1}^r P(-\delta_i)$ inherits a
$\mathbb{Z}^m$-grading from~$P$ in the natural way.
Again we say that~$F$ is
{\sem graded by}~$W$.

In order to be able to use these gradings in our algorithms,
we need some positivity assumptions.

\begin{definition}
Let~$P$ be graded by~$W$, and let $w_1,\dots,w_m$ be the rows
of~$W$.

\begin{items}
\item The grading given by~$W$ is called {\sem weakly positive}
if there exist integers $a_1,\dots,a_m$ such that
$a_1w_1 +\cdots +a_mw_m$ has all entries strictly positive.

\item The grading given by~$W$ is called {\sem positive} if
$\rk(W)=m$, if no column of~$W$ is zero, and
if the first non-zero entry in each column of~$W$ is positive.

\end{items}
\end{definition}

\begin{prop}
\label{gradedNAK}
Let~$P$ be weakly positively graded by~$W\!$, and let~$M$ be a
finitely generated graded $P$-module.

\begin{items}
\item We have $P_{W,0}=K$ and $\dim_K(M_{W,d})<\infty$
for every $d\in \mathbb{Z}^m$.

\item The graded version of Nakayama's lemma holds:
homogeneous elements $v_1,\dots,v_s\in M$ generate the module~$M$
if and only if their residue classes $\overline{v}_1,\dots,
\overline{v}_s$ generate the $K$-vector space $M/(x_1,\dots,x_n)M$.
In particular, every homogeneous system of generators of~$M$
contains a minimal one, and all irredundant homogeneous systems
of generators of~$M$ have the same number of elements  which is 
denoted by $\mu(M)$.

\end{items}
\end{prop}

The proof of this proposition uses standard computer algebra 
methods and is
contained in~\citep{KR2003}. For practical computations we need the
somewhat stronger notion of a positive grading. The usefulness of
positive gradings is illustrated by the following
characterizations. Recall that a module ordering~$\sigma$ on the
set of terms $\Tnr$ of the graded free module~$F$ is called 
{\bf degree compatible} or {\bf compatible
with} $\deg_W$ if the inequality $\deg_W(te_i)
>_{\mathtt{Lex}} \deg_W(t'e_j)$ implies $te_i >_\sigma t'e_j$ for all
$t,t'\in \mathbb{T}^n$ and all $i,j\in\{1,\dots,r\}$.

\begin{prop}
\label{positive=well}
Let~$P$ be graded by~$W$, where~$W$ has $\mathbb{Z}$-linearly
independent rows and non-zero columns. Then the following
conditions are equivalent.

\begin{items}
\item The grading on~$P$ given by~$W$ is positive.

\item The restriction of~$\mathtt{Lex}$ to the monoid
$\Gamma=\{d\in \mathbb{Z}^m \mid P_{W,d}\ne 0\}$
is a well-ordering, i.e.\ every non-empty subset
of~$\Gamma$ has a minimal element with respect
to~$\mathtt{Lex}$.

\item The restriction of~$\mathtt{Lex}$ to the monoid
$\Gamma=\{d\in \mathbb{Z}^m \mid P_{W,d}\ne 0\}$
is a term ordering, i.e.\ every element $d\in\Gamma$
satisfies $d >_{\mathtt{Lex}} 0$.

\item There exists a term ordering~$\tau$ on~$\mathbb{T}^n$
which is compatible with~$\deg_W$.

\item There exists a module term ordering~$\sigma$
on~$\Tnr$ which is compatible with the grading given
by~$W$.

\end{items}
\end{prop}

Again we refer to~\citep{KR2003} for a proof of this proposition. As
a consequence, it follows that positive gradings are weakly
positive. Moreover, in a positively graded setting, we can prove
the finiteness of various algorithms in the usual way, i.e.\ by
using the fact that there is no infinite sequence of homogeneous
elements of strictly decreasing degrees.

In the remaining part of this section, we use truncated
Gr\"obner bases to prove two very important technical tools,
namely Corollary~\ref{truncindegd} and
Corollary~\ref{addonegen}. We shall from now on assume that~$P$
is positively graded by $W\in\Mat_{m,n}(\mathbb{Z})$.
Moreover, we let $\delta_1,\dots,\delta_r\in\mathbb{Z}^m$,
we let~$M$ be a finitely generated graded submodule of the
graded free $P$-module $F=\dsum_{i=1}^r P(-\delta_i)$, and
we let~$\sigma$ be a module term ordering on~$\Tnr$,
the set of terms in~$F$.

The following notation will turn out to be convenient.
Given a subset~$S$ of a graded $P$-module and 
$d\in \mathbb{Z}^m$, we let $S_{\le d}=
\{v\in S\mid v\hbox{\ \rm homogeneous}$, $\deg_W(v)
\le_{\mathtt{Lex}} d\}$ and 
$S_d = \{v\in S\mid v$ homogeneous, $\deg_W(v) = d\}$.

\begin{definition}
Assume that $G=\{g_1,\dots, g_s\}$ is a homogeneous
$\sigma$-Gr\"obner basis of~$M\!$, and let $d\in \mathbb{Z}^m$.
Then the set $G_{\le d}$ is called a {\sem  $d$-truncated
Gr\"obner basis} of~$M\!$, or a Gr\"obner basis of~$M$
which has been {\sem truncated in degree~$d$}.
\end{definition}

For truncated Gr\"obner bases, we now prove a characterization
which is analogous to the Buchberger criterion in the usual
case. To this end, we need to explain what we mean by 
critical pairs and critical syzygies.

Given homogeneous elements $g_1,\dots,g_s\in M\setminus \{0\}$, we
let $d_i=\deg_W(g_i)$ for $i=1,\dots,s$, and we let~$F'$
be the graded free $P$-module $\dsum_{i=1}^s P(-d_i)$.
The canonical basis of~$F'$ will be denoted by $\{\epsilon_1,\dots,
\epsilon_s\}$. Notice that we have $\deg_W(\epsilon_i)=d_i$ for $i=1,\dots,s$.
Moreover, we write $\LM_\sigma(g_i) = c_i t_i e_{\gamma_i}$,
where $c_i\in K\setminus \{0\}$, where $t_i\in\mathbb{T}^n$, and where
$\gamma_i\in\{1,\dots,r\}$.

\begin{definition} A pair $(i,j)\in \{1,\dots,s\}$ such that
$1\le i<j\le s$ and $\gamma_i=\gamma_j$ is called a
{\sem critical pair} of $(g_1,\dots,g_s)$. The set of all critical
pairs of $(g_,\dots,g_s)$ is denoted by~$\mathbb{B}$.
For every critical pair $(i,j)\in\mathbb{B}$, the element
$\sigma_{ij}={\lcm(t_i,t_j) \over c_i t_i}\, \epsilon_i
- {\lcm(t_i,t_j) \over c_jt_j}\, \epsilon_j$ is a syzygy
of the pair $(\LM_\sigma(g_i),\LM_\sigma(g_j))$. It is called
the {\sem critical syzygy} associated to the critical pair
$(i,j)$. The set of all critical syzygies is denoted by~$\Sigma$.
\end{definition}

Clearly, a critical syzygy $\sigma_{ij}$ is a homogeneous element
of~$F'$ whose degree is precisely $\deg_W(\sigma_{ij})=
\deg_W(\lcm(t_i,t_j))+ \delta_{\gamma_i}$.
This degree equals the degree of the corresponding
\hbox{S-vector} $S_{ij}={\lcm(t_i,t_j) \over c_i t_i}\, g_i
- {\lcm(t_i,t_j) \over c_jt_j}\, g_j$ in~$F$.

For every critical pair $(i,j)\in\mathbb{B}$, we call $\deg_W(\sigma_{ij})$
the {\sem degree of the critical pair}. Then it makes sense to consider 
the set $\mathbb{B}_{\le d}$ for every given $d\in \mathbb{Z}^m\!$, 
and we observe that $\deg_W(\sigma_{ij}) \ge_{\mathtt{Lex}} \max\{d_i,d_j\}$
for all $(i,j)\in \mathbb{B}$.  Finally, we remind the reader that
$\NR_{\sigma,\Cal{G}}(v)$ denotes normal remainder, i.e.\ the result of
the division algorithm, as in \citep{KR2000}, Definition 1.6.7.
At this point, we are ready to formulate 
and prove the following characterization of truncated Gr\"obner bases.

\begin{prop}{\sem (Characterization of Truncated Gr\"obner
Bases)}\\
\label{charoftrunc}%
Let~$P$ be positively graded by $W\in\Mat_{m,n}(\mathbb{Z})$,
let $G=\{g_1,\dots,g_s\}$ be a set of non-zero homogeneous
vectors which generates a graded submodule~$M$ of~$\dsum_{i=1}^r
P(-\delta_i)$, and let $d\in \mathbb{Z}^m$. Then the following conditions
are equivalent.

\begin{items}
\item The set $G_{\le d}$ is a $d$-truncated $\sigma$-Gr\"obner
basis of~$M\!$.

\item For every homogeneous element $v\in M_{\le d}\setminus
\{0\}$, we have the relation 
$\LT_\sigma(v)\in \langle\LT_\sigma(g)\mid g\in G_{\le d}\rangle$.

\item For all pairs $(i,j)\in\mathbb{B}_{\le d}$, we have
$\NR_{\sigma, \mathcal{G}_{\le d}}(S_{ij})= 0$, where
$\mathcal{G}_{\le d}$ is the tuple obtained
from~$\mathcal{G}=(g_1,\dots,g_s)$
by deleting the elements of degree greater than~$d$.

\end{items}
\end{prop}

\proof
Without loss of generality, we may assume that $G_{\le d}=
\{g_1,\dots,g_{s'}\}$ for some $s'\le s$. It is clear that~a)
implies both~b) and~c). Now we show that~b) implies~a).
By the assumption, we can find terms $t'_{s'+1},\dots,t'_{s''}$
of degree greater than~$d$ such that the set
$\{\LT_\sigma(g_1),\dots,\LT_\sigma(g_{s'})\} \;\cup\;
\{t'_{s'+1},\dots,t'_{s''}\}$
is a system of generators of~$\LT_\sigma(M)$.
We choose homogeneous elements $h_{s'+1},\dots,h_{s''}$ in~$M$
such that $\LT_\sigma(h_i)=t'_i$ for $i=s'+1,\dots,s''$.
Then the set $\{g_1,\dots,g_{s'},\allowbreak
h_{s'+1},\dots,h_{s''}\}$ is a homogeneous \hbox{$\sigma$-Gr\"obner} 
basis of~$M$ with truncation $G_{\le d}$.

It remains to prove that c) implies b).
Let $v\in M_{\le d}$ be a homogeneous non-zero element.
Since $\{g_1,\dots,g_{s'}\}$ generates $\langle M_{\le d}\rangle$,
we can represent~$v$ as \hbox{$v = \sum_{i=1}^{s'}f_ig_i$}, where~$f_i$
is homogeneous of degree $\deg_W(v)-\deg_W(g_i)
\le_{\mathtt{Lex}} d$. In order to prove $\LT_\sigma(v) \in
\langle \LT_\sigma(g_1),\dots,\LT_\sigma(g_{s'})\rangle$,
it is enough to proceed as in the proof of 
Proposition~2.3.12 of~\citep{KR2000}, replacing $\mathcal{G}$ 
by~$\mathcal{G}_{\le d}$.
\endproof

This characterization has several useful applications.

\begin{cor}\label{truncoftrunc}
Let $G=\{g_1,\dots, g_s\}$ be a homogeneous $\sigma$-Gr\"obner
basis of the module~$M$, and let $d\in \mathbb{Z}^m$. Then $G_{\le d}$ is a
$d$-truncated $\sigma$-Gr\"obner basis of the
module $\langle M_{\le d}\rangle$.
\end{cor}

\proof
Since~$G$ is a set of generators of~$M$, the set~$G_{\le d}$
generates the module $\langle M_{\le d}\rangle$. From
Buchberger's Criterion we know that $\NR_{\sigma,\mathcal{G}}(S_{ij}) = 0$,
for all pairs $(i,j)\in \mathbb{B}$. If we have $\deg_W(S_{ij})
\le_{\tt Lex} d $ here, the elements of~$G$ involved in
the reduction steps $S_{ij} \TTo{G} 0$ all have degrees less than
or equal to~$d$. Hence we see that $\NR_{\sigma,\Cal{G}_{\le d}}
(S_{ij}) = 0$, and the proposition yields the claim.
\endproof

\begin{cor}\label{truncindegd}
Let $d\in\mathbb{Z}^m\!$, let the elements of the tuple
$\Cal{G}=(g_1,\dots, g_s)$ form a
$d$-truncated $\sigma$-Gr\"obner basis of~$M\!$,
and let $g_{s+1} \in F$ be a
homo\-geneous element of degree~$d$
such that $\LT_\sigma(g_{s+1})\notin  \langle \LT_\sigma(g_1),
\dots,\LT_\sigma(g_s)\rangle$. Then $\{g_1,\dots,g_{s+1}\}$
is a $d$-truncated Gr\"obner basis of~$M+\langle g_{s+1}\rangle$.
\end{cor}

\proof
In order to prove the claim, we check condition~c) of
the proposition. For $1\le i<j\le s$
such that $\deg_W(S_{ij})\le_{\tt Lex} d$, we have
$\NR_{\sigma,\Cal{G}}(S_{ij})= 0$ by the assumption and
by Proposition~\ref{charoftrunc}. For $i\in \{1,\dots,s\}$
such that $\deg_W(S_{i\,s{+}1})=~d$, the fact that
the pair $(i,s+1)$ has degree~$d$ implies that
$\LT_\sigma(g_{s+1})$ is a multiple of $\LT_\sigma(g_i)$,
in contradiction to the hypothesis.
\endproof

In the last part of this section, we prove an analogue of the
preceding corollary for minimal generators. Recall that
Proposition~\ref{gradedNAK}.b guarantees that all minimal systems
of generators have the same length in the positively graded
situation.

\begin{prop}\label{orderminimal}
Let $P$ be positively graded by~$W\in\Mat_{m,n}(\mathbb{Z})$, 
let~$M$ be a graded $P$-module generated 
by homogeneous elements $\{g_1,\dots,g_s\}$, and assume that
$\deg_W(g_1) \le_{\tt Lex} \deg_W(g_2) \le_{\tt Lex} \cdots 
\le_{\tt Lex} \deg_W(g_s)$.

\begin{items}
\item The set $\{g_1,\dots,g_s\}$ is a minimal system of generators
of~$M$ if and only if we have $g_i \notin \langle g_1,\dots,
g_{i-1}\rangle$ for $i=1,\dots,s$.

\item The set $\{g_i \mid i\in\{1,\dots,s\},\, g_i \notin
\langle g_1,\dots,g_{i-1}\rangle\}$ is a minimal system
of generators of~$M$.

\end{items}
\end{prop}

\proof
First we prove~a).
If $\{g_1,\dots, g_s\}$ is a minimal set of generators of~$M$,
then no relation of type $g_i \in \langle g_1,\dots, g_{i-1}\rangle$
holds, since otherwise we would have
$M = \langle g_1,\dots, g_{i-1}, g_{i+1}, \dots, g_s \rangle$.
Conversely, if $\{g_1,\dots, g_s\}$ is not a minimal set of
generators of~$M$, then there exists an index $i\in\{1,\dots,s\}$
such that $g_i\in \langle g_1,\dots, g_{i-1}, g_{i+1},\dots,
g_s\rangle$. Using Corollary~1.7.11 of~\citep{KR2000},
we obtain a representation
$g_i = \sum_{j\ne i}f_jg_j$, where $f_j\in P$ is homogeneous
of degree~$\deg_W(g_i)-\deg_W(g_j)$ for $j\in\{1,
\dots,i-1,i+1,\dots,s\}$.

Since $\deg_W(f_j)\ge_{\tt Lex}0$ for $f_j\ne 0$, 
we see that $\deg_W(g_i) <_{\tt Lex} \deg_W(g_j)$ implies $f_j=0$. 
Thus there are two possibilities. 
Either we have $\deg_W(g_i) >_{\tt Lex} \deg_W(g_j)$
for all~$j$ such that $f_j\ne 0$ or there exist
some indices~$j$ such that $\deg_W(g_j)=\deg_W(g_i)$.
In the first case, those indices~$j$
satisfy $j<i$ by the assumption that the multidegrees of
$g_1,\dots,g_s$ are ordered increasingly, and therefore we get
$g_i\in \langle g_1,\dots,g_{i-1}\rangle$. 
In the second case, the~$f_j$ corresponding 
to those indices~$j$ are in $K\setminus \{0\}$.
Let $j_{\max}= \max\{j\in \{1,\dots,s\}\mid
f_j \in K\setminus \{0\}\}$. We get the relation $g_{j_{\max}} \in
\langle g_1,\dots,g_{j_{\max} -1}\rangle$.
In both cases, we arrive at a contradiction to our hypothesis.

Now let us show~b). 
The set $S=\{g_i \mid i\in
\{1,\dots,s\},\, g_i \notin \langle g_1,\dots,g_{i-1}\rangle\}$
is a system of generators of~$M$, because an element~$g_i$
such that $g_i\in\langle g_1,\dots,g_{i-1}\rangle$ is also
contained in $\langle g_j\in S \mid 1\le j\le i-1\rangle$.
The fact that this system of generators is minimal
follows from~a).
\endproof

The following version is an immediate consequence of
part~a) of the proposition.

\begin{cor}\label{addonegen}
Let $N$ be a graded $P$-module, let $M$ be a submodule of~$N\!$,
let $\{g_1,\dots,g_s\}$ be a
minimal homogeneous system of generators of~$M$, and let
$g_{s+1}\in N\setminus M$ be a homogeneous
vector whose degree satisfies the inequality 
$\deg_W(g_{s+1})\ge_{\tt Lex} \max \{\deg_W(g_i)
\mid i=1,\dots,s\}$.
Then $\{g_1,\dots,g_{s+1}\}$ is a minimal system of generators
of the module $M+\langle g_{s+1}\rangle$.
In particular, we have $\mu(M+\langle g_{s+1}\rangle) =\mu(M)+1$.
\end{cor}


\section{Minimal Generators in a Reduced Gr\"obner Basis}
\label{Minimal Generators in a Reduced Groebner Basis}


From here on we use the following assumptions.
Let~$K$ be a field, and let $P=K[x_1,\dots,x_n]$ be a
polynomial ring over~$K$ which is positively graded
by a matrix~$W\in\Mat_{m,n}(\mathbb{Z})$. Then let $r\ge 1$,
let $\delta_1,\dots,\delta_r\in\mathbb{Z}^m$, and let~$M$
be a graded submodule of~$F=\dsum_{i=1}^r P(-\delta_i)$
which is generated by a set of non-zero homogeneous
vectors $\{v_1,\dots,v_s\}$. Furthermore, we choose
a module term ordering~$\sigma$ on the monomodule
of terms $\mathbb{T}^n\langle e_1,\dots,e_r\rangle$
in~$F$, and we let
$\Cal{V}=(v_1,\dots,v_s)$.

Our first goal is to describe an algorithm which computes a
homogeneous \hbox{$\sigma$-Gr\"obner} basis of~$M$ degree-by-degree
and a variant of this algorithm which also yields a minimal
system of generators of~$M$ contained in~$\Cal{V}$. This
part is classical and more or less ``well-known''. Then we
make good use of it in Theorem~\ref{BAMforRGB} for
minimalizing reduced Gr\"obner bases.

To ease the notation, we shall use the following convention:
whenever a vector~$g_i$ appears, we write
$\LM_\sigma(g_i)=c_it_ie_{\gamma_i}$, where
$c_i\in K\setminus \{0\}$, where $t_i\in\mathbb{T}^n$,
and where $\gamma_i\in\{1,\dots,r\}$.
For two indices $i,j$ such that $\gamma_i=\gamma_j$,
we let $\sigma_{ij}=\frac{\lcm(t_i,t_j)}{c_it_i}\,
\epsilon_i - \frac{\lcm(t_i,t_j)}{c_jt_j}\, \epsilon_j$ and $S_{ij}=
\frac{\lcm(t_i,t_j)}{c_it_i}\, g_i -
\frac{\lcm(t_i,t_j)}{c_jt_j}\, g_j$.

\begin{thm}{\bf (The Homogeneous Buchberger
Algorithm)}\\
\label{homBA}%
In the above situation, consider the following
instructions.

\begin{items}
\item[1)] Let $B=\emptyset$, $\Cal{W}=\Cal{V}$,
$\Cal{G}=\emptyset$, and let $s'=0$.

\item[2)] Let~$d$ be the smallest degree
with respect to~{\tt Lex} of an element of~$B$ or of~$\Cal{W}$.
Form~$B_d$ and~$\Cal{W}_d$, and delete their entries
from~$B$ and~$\Cal{W}$, respectively.

\item[3)] If~$B_d=\emptyset$, continue with step~6). Otherwise, chose a
pair $(i,j)\in B_d$ and remove it from~$B_d$.

\item[4)] Compute the S-vector~$S_{ij}$ and its normal remainder
$S'_{ij}=\NR_{\sigma,\Cal{G}}(S_{ij})$. If~$S'_{ij}=0$, continue with
step~3).

\item[5)] Increase $s'$ by one, append $g_{s'}=S'_{ij}$ to the
tuple~$\Cal{G}$, and append the set $\{(i,s')\mid 1\le i<s',\,\gamma_i=
\gamma_{s'}\}$ to the set~$B$. Continue with step~3).

\item[6)] If~$\Cal{W}_d=\emptyset$, continue with step~9). Otherwise,
choose a vector $v\in\Cal{W}_d$ and remove it from~$\Cal{W}_d$.

\item[7)] Compute $v'=\NR_{\sigma,\Cal{G}}(v)$. If~$v'=0$,
continue with step~6).

\item[8)] Increase~$s'$ by one, append $g_{s'}=v'$ to the
tuple~$\Cal{G}$, and append the set $\{(i,s')\mid 1\le i<s',\,
\gamma_i=\gamma_{s'}\}$ to the set~$B$. Continue with
step~6).

\item[9)] If $B=\emptyset$ and $\Cal{W}=\emptyset$, return the
tuple~$\Cal{G}$ and stop. Otherwise, continue with step~2).
\end{items}

\noindent This is an algorithm which returns a 
$\sigma$-Gr\"ob\-ner basis $\Cal{G}$ of~$M$,
where the tuple~$\Cal{G}$ consists of homogeneous 
vectors having non-decreasing multidegrees.
\end{thm}

The proof of this theorem is standard Computer Algebra and
is for instance contained in~\citep{KR2003}.

\begin{remark}\label{comptruncGB}
Let us add some observations about this algorithm.

\begin{items}
\item If we interrupt its execution after some degree~$d_0$ is
finished, the tuple~$\Cal{G}$ is a $d_0$-truncated
Gr\"obner basis of~$M$. Consequently, we can compute
truncated Gr\"obner bases efficiently.
Moreover, in this case it suffices to append only the pairs
$\{(i,s')\mid 1\le i< s',\, \gamma_i=\gamma_{s'},\,
\deg_W(\sigma_{is'})\le_{\tt Lex} d_0\}$ to the set~$B$
in steps~5) and~8). The reason is that pairs of higher
degree are never processed anyway, since we stop the computation
after finishing degree~$d_0$.

\item It is not required that~$\sigma$
is a degree compatible module term ordering. The reason is that,
during the computation of the Gr\"obner basis, only comparisons
of terms in the support of a homogeneous vector are performed.
Thus these terms have the same degree, and it does not matter
whether~$\sigma$ is degree compatible or not.

\item The Homogeneous Buchberger Algorithm can also be viewed
as a special version of the usual Buchberger Algorithm where
we use a suitable selection strategy.

\end{items}
\end{remark}

The following variant of the Homogeneous Buchberger Algorithm
computes a minimal system of generators of~$M$ contained in the
given set of generators while computing a Gr\"obner basis. It provides an
efficient method for finding minimal systems of generators.

\medbreak\goodbreak

\begin{cor}{\sem (Buchberger Algorithm with Minimalization)}\\
\label{BAM}\nobreak
In the situation of the theorem, consider the following instructions.

\begin{items}
\item[1')]Let $B=\emptyset$, $\Cal{W}=\Cal{V}$,
$\Cal{G}=\emptyset$, $s'=0$, and $\Cal{V}_{\min}=\emptyset$.

\item[2)] Let~$d$ be the smallest degree
with respect to~{\tt Lex} of an element of~$B$ or of~$\Cal{W}$.
Form~$B_d$ and~$\Cal{W}_d$, and delete their entries
from~$B$ and~$\Cal{W}$, respectively.

\item[3)] If~$B_d=\emptyset$, continue with step~6). Otherwise, chose a
pair $(i,j)\in B_d$ and remove it from~$B_d$.

\item[4)] Compute the S-vector~$S_{ij}$ and its normal remainder
$S'_{ij}=\NR_{\sigma,\Cal{G}}(S_{ij})$. If~$S'_{ij}=0$, continue with
step~3).

\item[5)] Increase $s'$ by one, append $g_{s'}=S'_{ij}$ to the
tuple~$\Cal{G}$, and append the set $\{(i,s')\mid 1\le i<s',\,\gamma_i=
\gamma_{s'}\}$ to the set~$B$. Continue with step~3).

\item[6)] If~$\Cal{W}_d=\emptyset$, continue with step~9). Otherwise,
choose a vector $v\in\Cal{W}_d$ and remove it from~$\Cal{W}_d$.

\item[7)] Compute $v'=\NR_{\sigma,\Cal{G}}(v)$. If~$v'=0$,
continue with step~6).

\item[8')] Increase~$s'$ by one, append $g_{s'}=v'$ to the
tuple~$\Cal{G}$, append~$v$ to the tuple~$\Cal{V}_{\min}$, and
append $\{(i,s')\mid 1\le i<s',\,
\gamma_i=\gamma_{s'}\}$ to the set~$B$. Continue with step~6).

\item[9')] If $B=\emptyset$ and $\Cal{W}=\emptyset$, return the
pair~$(\Cal{G},\Cal{V}_{\min})$ and stop. Otherwise, continue with
step~2).

\end{items}

\noindent This is an algorithm which returns a 
pair~$(\Cal{G},\Cal{V}_{\min})$ such
that~$\Cal{G}$ is a tuple of homogeneous vectors which are a
$\sigma$-Gr\"obner basis of~$M$, and~$\Cal{V}_{\min}$ is a subtuple
of~$\Cal{V}$ of homogeneous vectors which are a minimal system of
generators of~$M$.
\end{cor}

\proof
In view of the theorem, we only have to show that the elements
in~$\Cal{V}_{\min}$ are a minimal set of generators of~$M$. Since
the algorithm is finite, it operates in only finitely many
degrees~$d$. Therefore it suffices to prove by induction on~$d$
that~$\Cal{V}_{\min}$ contains a minimal system of generators
of~$\langle M_{\le d}\rangle$ after the algorithm has finished
working on elements of degree~$d$.

This is clearly the case at the outset. Suppose it is true for the
last degree treated before~$d$. Inductively, we can show that the
elements of~$\Cal{G}$ continue to be contained in the
module~$\langle M_{<d} \rangle$ while we are looping through
steps~3), 4), and~5) of the algorithm. Namely, every time an
element of the form $\NF_{\sigma,\Cal{G}}(S_{ij})$ is added
to~$\Cal{G}$, it is clearly contained in the module generated by
the previous elements of~$\Cal{G}$. Furthermore, by part~a) of the 
remark following Theorem~\ref{homBA}, 
the elements of the tuple~$\Cal{G}$
form a $d$-truncated Gr\"obner basis of $\langle M_{<d}\rangle$
after we have finished looping through steps~3), 4), and~5), i.e.\
when we have treated all pairs of degree~$d$.

Now let $\Cal{W}_d=(w_1,\dots,w_\ell)$, and let the numbering of
these vectors correspond to the order in which they are chosen in
step~6). We show that, for each application of steps~6), 7),
and~8'), the elements of~$\Cal{V}_{\min}$ continue to be a minimal
system of generators of the module they generate, and that this
module always agrees with the one generated by the elements
of~$\Cal{G}$. Furthermore, the elements of~$\Cal{G}$ are always a
$d$-truncated $\sigma$-Gr\"obner basis of that module.

When a new vector~$v=w_i$ is chosen in step~6), there are two
possibilities. If~$v'=0$ in step~7), then~$v$ is already contained
in the module~$M'$ generated by the elements of~$\Cal{V}_{\min}$.
Otherwise, the vector~$v'$ is not contained in~$M'\!$, since the
elements of~$\Cal{G}$ are a $d$-truncated $\sigma$-Gr\"obner basis
and we can apply the Submodule Membership Test (see~\citep{KR2000},
Proposition 2.4.10.a). In that case, the elements of~$\Cal{V}_{\min}$, together
with~$v$, form a minimal system of generators of the module
$M'+\langle v\rangle = M'+\langle v'\rangle$ by
Corollary~\ref{addonegen}. Moreover, the elements of~$\Cal{G}$,
together with~$v'$, form a $d$-truncated $\sigma$-Gr\"obner basis
of~$M'+\langle v'\rangle$ by Corollary~\ref{truncindegd}.

Altogether, it follows that, after degree~$d$ is finished, the
elements of~$\Cal{V}_{\min}$ are a minimal system of
generators of~$\langle M_{\le d}\rangle$, as we wanted to show.
\endproof

\begin{remark}
Let us collect some observations about this
algorithm.

\begin{items}
\item If we are only interested in a minimal system of
generators of~$M$ (and not in a Gr\"obner basis), we can stop the
algorithm after we have completed degree
$d_{\max}=\max\{\deg(v_i)\mid 1\le i\le s\}$. In this case it
suffices to append only the pairs $\{(i,s') \mid 1\le i< s',\,
\gamma_i=\gamma_{s'},\,
\deg_W(\sigma_{is'})\le_{\tt Lex} d_{\max}\}$ to the set~$B$
in steps~5) and~8').

\item
In addition, we could alter step~8') and append the vector~$v'$
instead of~$v$ to the list~$\Cal{V}_{\min}$. Then~$\Cal{V}_{\min}$
would still contain a minimal homogeneous set of generators of~$M$
when the computation ends. These generators would not be contained
in the initial tuple~$\Cal{V}$ anymore, but they would have the
additional property that each vector is fully reduced against the
previous ones.
\end{items}
\end{remark}

The final part of the section is devoted to a result 
which will be essential for our discussion of the minimalization 
of the critical pairs. 
Namely, we are going to apply the algorithm of Corollary~\ref{BAM} 
to a reduced Gr\"obner basis and improve it significantly 
in that case. The main differences between both algorithms occur
in step~7), where it suffices to
compare terms instead of computing normal remainders, and in step~8),
where we append~$v$ to both~$\Cal{G}$ and~$\Cal{V}_{\min}$.

\begin{thm}{\sem (Minimal Generators in a Reduced
Gr\"obner Basis)}\\
\label{BAMforRGB}%
In the situation of Theorem~\ref{homBA},
let~$\Cal{V}=(v_1,\dots,v_s)$
be the reduced $\sigma$-Gr\"obner basis of~$M$.
Consider the following instructions.

\begin{items}
\item[1)] Let $B=\emptyset$, $\Cal{W}=\Cal{V}$,
$\Cal{G}=\emptyset$, $s'=0$, and $\Cal{V}_{\min}=\emptyset$.

\item[2)] Let~$d$ be the smallest degree
with respect to~{\tt Lex} of an element
of~$B$ or of~$\Cal{W}$.
Form~$B_d$ and~$\Cal{W}_d$, and delete their entries
from~$B$ and~$\Cal{W}$, respectively.

\item[3)] If $B_d=\emptyset$, continue with step 6).
Otherwise, choose a pair $(i,j)\in B_d$ and remove it
from~$B_d$.

\item[4)] Compute $S'_{ij}=\NR_{\sigma,\Cal{G}}(S_{ij})$.
If $S'_{ij}=0$, continue with step~3).

\item[5)] Increase $s'$ by one, append
$g_{s'}=S'_{ij}$ to the tuple~$\Cal{G}$, append the following  
set $\{(i,s')\mid 1\le
i<s',\, \gamma_i=\gamma_{s'}\}$ to~$B$, and continue with
step~3).

\item[6)] If~$\Cal{W}_d=\emptyset$, continue with step~9).
Otherwise, choose $v\in \Cal{W}_d$ and remove it
from~$\Cal{W}_d$.

\item[7)] If~$\LT_\sigma(v)=\LT_\sigma(g)$ for some
$g\in\Cal{G}$, then replace the element~$g$ in~$\Cal{G}$
by~$v$. Continue with step~6).

\item[8)] Increase $s'$ by one, append~$g_{s'}=v$
to the tuples~$\Cal{G}$ and~$\Cal{V}_{\min}$, and append $\{(i,s')\mid 1\le
i<s',\, \gamma_i=\gamma_{s'}\}$ to the set~$B$. Continue with
step~6).

\item[9)] If $B=\emptyset$ and $\Cal{W}=\emptyset$,
return~$\Cal{V}_{\min}$ and stop. Otherwise, continue with step~2).

\end{items}

\noindent This is an algorithm which computes
a subtuple~$\Cal{V}_{\min}$ of~$\Cal{V}$ such
that~$\Cal{V}_{\min}$ is a minimal system of generators of~$M$.
\end{thm}

\proof
It suffices to show that this procedure has the same effect
as running the algorithm of Corollary~\ref{BAM} on~$\Cal{V}$.

First we use induction on~$d$ to show that, after we have
finished some degree~$d$, the tuple~$\Cal{G}$ has the same
elements as~$\Cal{V}_{\le d}$. Every element of~$\Cal{V}_d$
is appended to~$\Cal{G}$ at some point in step~7) or~8). On
the other hand, if an element~$g_{s'}$ is put into~$\Cal{G}$
in step~5), it has a leading term which is not a multiple of
an element of~$\Cal{V}_{<d}$. Hence it is swapped out
of~$\Cal{G}$ at some point in step~7).

Next we note that, after we have finished cycling through
steps~3), 4), and~5) in degree~$d$, the tuple~$\Cal{G}$ is a
$d$-truncated minimal $\sigma$-Gr\"obner basis of~$M_{<d}$.

Now we turn our attention to the loop described in steps~6),
7) and~8). Notice that the effect of steps~7) and~8) is
independent of the order in which we choose the elements
$v\in\Cal{W}_d$ in step~6). Hence we can assume for the
purposes of this proof that we always choose the vector~$v$
in~$\Cal{W}_d$ which has the minimal leading term with
respect to~$\sigma$. With this assumption, we show
inductively that when we run steps~7) and~8) for some
element $v\in \Cal{W}_d$, at each point the elements
in~$\Cal{G}$ are a minimal $\sigma$-Gr\"obner basis of the
module they generate, and the elements of~$\Cal{V}_{\min}$
are a minimal system of generators of that module.

For the induction step, we have to consider two cases:
either~$v$ is swapped into~$\Cal{G}$ in step~7) or appended
to both~$\Cal{G}$ and~$\Cal{V}_{\min}$ in step~8). In the first
case, it suffices to show that the module generated by
the elements of~$\Cal{G}$ does
not change when we perform the swap, i.e.\ that the
difference $v-g$ is contained in this module.
This follows from the observations that $\LT_\sigma(v-g)
<_\sigma \LT_\sigma(v)$ and all elements~$\tilde v$
in~$\Cal{V}$ such  that $\LT_\sigma(\tilde v) <_\sigma
\LT_\sigma(v)$ are already in~$\Cal{G}$. Since $v-g
\TTo{\Cal{V}}0$, we have $v-g\TTo{\Cal{G}}0$.
In the second case, it is clear that~$\Cal{G}$ continues to
be a minimal Gr\"obner basis of the module it generates by
Corollary~\ref{truncindegd}, and~$\Cal{V}_{\min}$ continues to be a
minimal system of generators of that module by
Corollary~\ref{addonegen}.

Finally, we note that in step~8) we can append~$v$
to~$\Cal{G}$ without passing to the normal remainder,
since~$v$ is an element of a reduced Gr\"obner basis and
thus irreducible.
\endproof

\begin{remark}\label{remMRGB}
Let us make some observations about the
preceding algorithm.

\begin{items}
\item The proof of the proposition shows that the algorithm
reconstructs the given reduced Gr\"obner basis
inside~$\Cal{G}$, and that~$\Cal{G}_{\le d}$ has the same
elements as~$\Cal{V}_{\le d}$ after some degree~$d$ is
finished.

\item Moreover, we note that in step~4) it is not necessary to
compute the normal remainder $\NR_{\sigma,\Cal{G}}(S_{ij})$.
Rather, it suffices to perform a full leading term
reduction.

\item The different elements $\NR_{\sigma,
\Cal{G}}(S_{ij})$ computed in step~4) and the elements
$v \in \Cal{V}_d$ which are swapped into $\Cal{G}$ by
step~7) are in $1-1$ correspondence, since every new element
computed in step~4) must have a new leading term in the
leading term module of~$M$. This new leading term must be
the leading term of an element in the reduced Gr\"obner
basis, hence it is swapped.

\end{items}
\end{remark}

\newpage 

\section{Minimalizing the Critical Syzygies}
\label{Minimalizing the Critical Syzygies}

In this section we continue to use the assumptions
and notation of the previous section. If we look at
Theorem~\ref{homBA} and
its proof, we can see that instead of treating all
pairs $(i,j)$ such that~$\sigma_{ij}$ is contained
in the set of critical syzygies~$\Sigma$, it would be enough to treat
those pairs corresponding to a subset
$\Theta\subseteq\Sigma$ which is a minimal system of
generators of~$\Syz_P(c_1 t_1 e_{\gamma_1},\dots,
c_s t_s e_{\gamma_s})$.

In order to find~$\Theta$, we observe that the application
of two of the rules for killing critical pairs
given in~\citep{GM1987} produces a minimal Gr\"obner basis
of the module~$\Syz_P(c_1 t_1 e_{\gamma_1},\dots,
c_s t_s e_{\gamma_s})$ contained in
the set~$\Sigma$. From this we derive the idea
to find~$\Theta$ by applying Theorem~\ref{BAMforRGB}.
We need the following definition.

\begin{definition} On the set of terms
$\mathbb{T}^n\langle \epsilon_1,\dots,\epsilon_s\rangle$ 
in~$\dsum_{i=1}^s P(-d_i)$ we define a relation~$\tau$ by letting
$$
t\, \epsilon_i \ge_\tau t'\, \epsilon_j \;\iff\;
\cases{t\,t_i\,e_{\gamma_i} >_\sigma
t'\,t_j\,e_{\gamma_j},\hbox{\quad \rm or} & \cr
t\,t_i\,e_{\gamma_i} =
t'\,t_j\,e_{\gamma_j} \hbox{\quad \rm and\quad }
i\ge j &}
$$
for $t,t'\in\mathbb{T}^n$ and $i,j\in\{1,\dots,s\}$. As
in~\citep{KR2000}, Lemma 3.1.2, it follows that~$\tau$ is a module term
ordering. It is called the term ordering {\sem induced} by the
tuple $(t_1 e_{\gamma_1},
\dots, t_s e_{\gamma_s})$ and by~$\sigma$.
\end{definition}

By~\citep{KR2000}, Proposition 3.1.3, the set~$\Sigma$ is a 
$\tau$-Gr\"obner basis of the module~$\Syz_P(c_1 t_1 e_{\gamma_1},
\dots, c_s t_s e_{\gamma_s})$. Moreover, $\sigma_{ij}$ is a 
homogeneous element of~$\dsum_{i=1}^s P(-d_i)$ of degree
$\deg_W(\sigma_{ij})=\deg(\lcm(t_i,t_j)) + \delta_{\gamma_i}$. 
For all $i,j\in\{1,\dots,s\}$, we let $t_{ij}={\lcm(t_i,t_j)\over t_i}$.
Now the main result of~\citet{GM1987} reads as follows.

\begin{prop}\label{GMrules}
Consider the following instructions.

\begin{rules}
\item[\sc Rule 1.] Delete in~$\Sigma$ all elements~$\sigma_{jk}$
such that there exists an index
$i$ in the set $\{1, \dots, j-1\}$ such that~$t_{ki}$ divides~$t_{kj}$. Call
the resulting set~$\Sigma'$.

\item[\sc Rule 2.] Delete in~$\Sigma'$ all
elements~$\sigma_{ik}$ such that there exists an index
$j$ in the set $\{i+1,\dots,k-1\}$ such that~$t_{kj}$ properly
divides~$t_{ki}$. Call the resulting set~$\Sigma''$.

\item[\sc Rule 3.] Delete in~$\Sigma''$ all
elements~$\sigma_{ij}$ such that there exists an index
$k$ in the set $\in\{j+1,\dots,s\}$ such that~$t_{ik}$ properly
divides~$t_{ij}$ and $t_{jk}$ properly divides~$t_{ji}$.
Call the resulting set~$\Sigma'''$.
\end{rules}

Then the set~$\Sigma'''$ still generates
$\Syz_P(c_1 t_1 e_{\gamma_1},\dots, c_s t_s e_{\gamma_s})$.
\end{prop}

\begin{remark} Let us interpret the previous proposition
in another way. For $1\le i<j\le s$ such that
$\gamma_i=\gamma_j$, we have $\LT_\tau(\sigma_{ij})=
t_{ji}\epsilon_j$. Hence Rules~1 and~2 can be restated
as follows.

\begin{rulesp}
\item[\sc Rule 1'.] Delete in~$\Sigma$ all
elements~$\sigma_{ij}$ such that there exists
an element~$\sigma_{i'j}$ such that $\LT_\tau
(\sigma_{ij})$ is a proper multiple
of~$\LT_\tau(\sigma_{i'j})$.

\item[\sc Rule 2'.] If, among the remaining elements,
there are elements~$\sigma_{ij}$,
$\sigma_{i'j}$ such that $\LT_\tau(\sigma_{ij})=
\LT_\tau(\sigma_{i'j})$, then delete the one having the
larger index $\max\{i,i'\}$.
\end{rulesp}
From Rules~1' and~2' it follows that
the set~$\Sigma''$ is a minimal $\tau$-Gr\"obner
basis of the module $\Syz_P(c_1 t_1 e_{\gamma_1},\dots,
c_s t_s e_{\gamma_s})$, i.e.\ the leading terms of the
elements of~$\Sigma''$ minimally generate the leading term 
module.
\end{remark}

In general, it is not true that $\Sigma''$ is a minimal system of
generators of the module $\Syz_P(c_1 t_1 e_{\gamma_1},\dots,
c_s t_s e_{\gamma_s})$, as our next example shows.
(For another example, see~\citep{GM1987}, 3.6.)

\begin{example}\label{ex1}
Let $P=\mathbb{Q}[x,y,z]$ be standard graded, 
let $r=1$, $s=4$ and $t_1=x^3z^2$,
$t_2=x^3y^4$, $t_3=y^5z^2$, $t_4=x^2y^5z$. Then we
get $\sigma_{12}=y^4 \epsilon_1-z^2 \epsilon_2$, $\sigma_{13}=
y^5 \epsilon_1 -x^3 \epsilon_3$, $\sigma_{14}= y^5 \epsilon_1 
-xz \epsilon_4$, $\sigma_{23}=yz^2 \epsilon_2- x^3 \epsilon_3$, 
$\sigma_{24}=yz \epsilon_2 - x \epsilon_4$, and 
$\sigma_{34}=x^2\epsilon_3 - z\epsilon_4$.
By applying Rules~1 and~2, we get the minimal
$\tau$-Gr\"obner basis $\Sigma''=\{ \sigma_{12},
\sigma_{24}, \sigma_{34}, \sigma_{13}\}$
of~$\Syz_P(t_1,t_2,t_3,t_4)$, since $\LT_\tau(\sigma_{23})
=\LT_\tau(\sigma_{13})$ and $\LT_\tau(\sigma_{14})=
z\cdot \LT_\tau(\sigma_{24})$.
Now we use Rule~3 and find $\Sigma'''=\Sigma''$,
but~$\Sigma'''$ is not a minimal system of generators
of~$\Syz_P(t_1,t_2,t_3,t_4)$, since we have 
$\sigma_{13}= y \sigma_{12} +z \sigma_{24}
-x \sigma_{34}$.
\end{example}

Before continuing, let us introduce a new notion. If we
have an element~$\sigma_{ij}$ and perform a reduction step
$\sigma_{ij}
\TTo{ct\sigma_{i'j}} c't' \epsilon_i + c''t'' \epsilon_{i'}$, 
where $c,c',c''\in~K$ and $t,t',t'' \in\mathbb{T}^n$, 
we call this a {\sem head reduction
step}. (Notice that the ~$j$-indices have to match!) Similarly, we
can define a {\sem tail reduction step} as follows: $\sigma_{ij}
\TTo{ct \sigma_{i'i}} c't' \epsilon_{i'} + c'' t'' \epsilon_j$.
It is clear that a tail reduction step does not change the
leading term of the element.

\begin{prop}\label{redSigma}
The set $\widetilde{\Sigma}=
\{ - c_j\cdot \sigma_{ij} \mid \sigma_{ij}
\in\Sigma''\}$ is the reduced $\tau$-Gr\"obner
basis of the module $\Syz_P(c_1 t_1 e_{\gamma_1},\dots,
c_s t_s e_{\gamma_s})$.
\end{prop}

\proof
Since passing from~$\Sigma''$ to $\widetilde{\Sigma}$ is equivalent
to normalizing the leading coefficients, and since~$\Sigma''$ is a
minimal $\tau$-Gr\"obner basis, it remains to show that no tail
reductions are  possible among the elements
of~$\mathstrut\widetilde{\Sigma}$. But if we perform a tail
reduction on some element of~$\mathstrut\widetilde{\Sigma}$, we get
an element of the form $\tilde c\,\tilde t\,\sigma_{i'j}$ such that
$i'<i$. Here we have to have $\tilde t=1$, since $\sigma_{ij}$ is
part of a minimal Gr\"obner basis. Now we obtain a contradiction to
the minimality of~$i$ in Rule~2'.
\endproof

\begin{remark}\label{remBAMapl}
Let us apply the algorithm of Theorem~\ref{BAMforRGB}
to the set~$\mathstrut\widetilde{\Sigma}$.
We make the following observations.

\begin{items}
\item A {\it pair of pairs}, i.e.\ a critical pair between
two elements $\sigma_{ij},\, \sigma_{i'j'}$ yields an
S-vector $S_{((i,j),(i',j'))}= c\,t\,\sigma_{ij} -
c'\,t'\, \sigma_{i'j'}$ such that $c,c'\in K$ and $t,t'\in
\mathbb{T}^n$ and $j=j'$, since the two leading terms have to cancel.
Without loss of generality, let $i<i'$. Then the result is $\tilde
c\,\tilde t\, \sigma_{ii'}$ for some $\tilde c\in K$ and $\tilde
t\in\mathbb{T}^n$.
The degree of such a pair of pairs is
 \begin{eqnarray*}
\deg_W(S_{((i,j),(i',j))}) & = &
\deg_W(\tilde t)+\deg_W(\sigma_{ii'})\\
& = &  \deg_W(\tfrac{\lcm(t_i,t_{i'},t_j)}{t_j}) + \deg_W(\epsilon_j)\\
& = & \deg_W(\lcm(t_i,t_{i'},t_j)) +\delta_{\gamma_j}
\end{eqnarray*}

\item During the course of the algorithm,
a new Gr\"obner basis element can only be obtained
from a pair of pairs if~$\tilde t=1$.
This is equivalent to $\gcd(t_{ij},t_{i'j})=1$.

\end{items}
\end{remark}

Now we are ready to optimize the minimalization of the critical
syzygies. To ease the notation, we shall minimalize the set~$\Sigma''$
instead of~$\widetilde{\Sigma}$. The lack of the normalization of
the leading coefficients is clearly of no consequence. We need the
following lemma.

\begin{lem}
Let $1\le i<j<m\le s$ and $i'\in \{1,\dots,j-1\}
\setminus \{i\}$. Suppose there are terms $t,t',t''\in
\mathbb{T}^n \setminus \{1\}$ such that $\sigma_{ij}=
\sigma_{ii'}+t\, \sigma_{i'j}=t'\, \sigma_{im} - t''\,
\sigma_{jm}$ and $\sigma_{i'm}=t\, \sigma_{i'j} + t''\,
\sigma_{jm}$. Then $t$, $t'$, and $t''$ are pairwise
coprime.

More precisely, given $\kappa\in \{1,\dots,n\}$, we
define~$\alpha=\deg_{x_\kappa}(t_i)$, $\alpha'=\deg_{x_\kappa}
(t_{i'})$, $\beta=\deg_{x_\kappa}(t_j)$,
and $\gamma=\deg_{x_\kappa}(t_m)$. Then one of the following
four cases occurs.

\begin{items}
\item[1)] We have $\alpha=\gamma > \beta$ and $\alpha > \alpha'$.

\item[2)] We have $\alpha' = \beta>\gamma$ and $\alpha'>\alpha$.

\item[3)] We have $\alpha=\alpha'>\beta$ and $\alpha>\gamma$.

\item[4)] We have $\alpha=\alpha'=\beta>\gamma$ or
$\alpha=\beta=\gamma>\alpha'$ or $\alpha'=\beta=\gamma>\alpha$.
\end{items}
\end{lem}

\proof
Comparing coefficients in the given equations yields the 
following equalities
$\lcm(t_i,t_j)=\lcm(t_i,t_{i'})=\lcm(t_{i'},t_m)=
t\,\lcm(t_{i'},t_j) = t'\,\lcm(t_i,t_m)= t''\,
\lcm(t_j,t_m)$. Thus the exponent of~$x_\kappa$ in these
terms satisfies $\max\{ \alpha,\beta\} = \max\{\alpha,
\alpha'\} = \max\{ \alpha',\gamma\} = \deg_{x_\kappa}
(t) + \max\{ \alpha',\beta\} = \deg_{x_\kappa}(t')+
\max\{ \alpha,\gamma\} = \deg_{x_\kappa}(t'')+
\max\{ \beta,\gamma\}$. We distinguish the following
four cases.

Case 1: Suppose that~$x_\kappa$ divides~$t$. In this case,
$\max\{ \alpha,\alpha'\} > \max\{ \alpha',
\beta\}$ yields $\alpha>\alpha'$ and $\alpha>\beta$.
Then $\alpha=\max\{\alpha,\alpha'\}=\max\{\alpha',\gamma\}$
shows $\alpha=\gamma$, i.e.\ we have the inequalities
stated in case~1) of the claim. Furthermore, it follows
that $\gamma=\max\{\alpha,\gamma\}=\max\{\beta,\gamma\}$,
i.e.\ that~$x_\kappa$ divides neither~$t'$ nor~$t''$.

Case 2: Suppose that~$x_\kappa$ divides~$t'$. In this
case, $\max\{\alpha,\alpha'\}>\max\{ \alpha,\gamma\}$
yields $\alpha'>\alpha$ and $\alpha'>\gamma$.
Then $\max\{\alpha,\beta\} = \max\{\alpha,\alpha'\}$
shows $\alpha'=\beta$, i.e.\ we have the inequalities
stated in case~2) of the claim. Furthermore, it follows
that $\beta= \max\{\alpha',\beta\} = \max\{\beta,\gamma\}$,
i.e.\ that~$x_\kappa$ divides neither~$t$ nor~$t''$.

Case 3: If~$x_\kappa$ divides~$t''$, we argue analogously
and obtain the inequalities stated in~3) as well as
the fact that~$x_\kappa$ divides neither~$t$ nor~$t'$.

Case 4: If~$x_\kappa$ divides neither~$t$ nor~$t'$
nor~$t''$, an easy case-by-case argument yields the
possibilities listed in~4).
\endproof

\begin{prop}{\sem (Minimalization of the
Critical Syzygies)}\\
\label{SigmaBAM}%
\noindent Let $\Sigma''$ be the $\tau$-Gr\"obner
basis of~$\Syz_P(c_1t_1e_{\gamma_1},\dots,c_st_se_{\gamma_s})$
defined in Proposition~\ref{GMrules}. 
Consider the following instructions.

\begin{items}
\item[1)] Let $\Cal{B}^\ast=\emptyset$, $\Cal{W}=\Sigma''$,
$\Cal{A}=\emptyset$, and $\Theta=\emptyset$.

\item[2)] For all $\sigma_{ij},\sigma_{i'j}\in\Sigma''$
such that $1\le i<i'<j\le s$,
form the S-vector $S_{((i,j),(i',j))}=
\tilde t\, \sigma_{ii'}$, where $\tilde t\in\mathbb{T}^n$.
If $\tilde t=1$, append~$\sigma_{ii'}$ to~$\Cal{B}^\ast$.

\item[3)] Let~$d$ be the smallest degree with respect
to~{\tt Lex} of an element of~$\Cal{B}^\ast$ or~$\Cal{W}$.
Form $\Cal{B}^\ast_d$ and~$\Cal{W}_d$, and delete
their entries from~$\Cal{B}^\ast$ and $\Cal{W}$, respectively.

\item[4)] If $\Cal{B}^\ast_d=\emptyset$, continue with step~11).
Otherwise, choose an element $\sigma_{ij}\in \Cal{B}^\ast_d$
and remove it from~$\Cal{B}^\ast_d$.

\item[5)] If $\LT_\tau(\sigma_{ij})\in\LT_\tau(\Cal{A}_d)$,
then continue with step~4).

\item[6)] If $\LT_\tau(\sigma_{ij}) = \LT_\tau(\sigma_{i'j})$
for some element $\sigma_{i'j}\in\Cal{W}_d$, then remove
$\sigma_{i'j}$ from~$\Cal{W}_d$, append it to~$\Cal{A}$, and
continue with step~4).

\item[7)] Find $\sigma_{i'j}\in \Cal{A}_{<d}$ such that
$t_{ji}$ is a multiple of~$t_{ji'}$.
Then perform the head reduction step
$\sigma_{ij} \TTo{\smash{\sigma_{i'\!j\,}}}
\tilde t \,\sigma_{k\ell}$, where $\tilde t\in\mathbb{T}^n$,
where $k=\min\{i,i'\}$, and where $\ell=\max\{i,i'\}$.
If $\tilde t\ne 1$, continue with step~4).

\item[8)] If $\LT_\tau(\sigma_{k\ell})\in \LT_\tau(\Cal{A}_d)$,
then continue with step~4).

\item[9)] If $\LT_\tau(\sigma_{k\ell})=
\LT_\tau(\sigma_{k'\ell})$ for some element
$\sigma_{k'\ell}\in \Cal{W}_d$,
then remove the element $\sigma_{k'\ell}$ from~$\Cal{W}_d$,
append~it to~$\Cal{A}$, and continue with step~4).

\item[10)] If $\sigma_{k\ell}\in \Cal{B}^\ast_d$,
then delete $\sigma_{k\ell}$ in~$\Cal{B}^\ast_d$
and continue with step~7), applied to this element.
Otherwise continue with step~4).

\item[11)] Append $\Cal{W}_d$ to~$\Cal{A}$ and to~$\Theta$.

\item[12)] If $\Cal{B}^\ast=\emptyset$ and $\Cal{W}=\emptyset$,
return~$\Theta$ and stop. Otherwise, continue with step~3).
\end{items}

\noindent This is an algorithm which computes
a subset~$\Theta\subseteq \Sigma''$ such that~$\Theta$ is a minimal
system of generators
of~$\Syz_P(c_1t_1e_{\gamma_1},\dots,c_st_se_{\gamma_s})$.
\end{prop}

\proof
It suffices to show that the given instructions define an
optimization of the application of Theorem~\ref{BAMforRGB} to the
set~$\Sigma''$. The tuple~$\Cal{A}$ corresponds to~$\Cal{G}$ there,
$\Theta$~corresponds to~$\Cal{V}_{\min}$, and $\Cal{B}^\ast$
corresponds to~$B$.

The first significant difference occurs in step~2).
Instead of producing the pairs of pairs inductively each
time we find a new Gr\"obner basis element, we precompute them
all at once. This is possible, since we know from
Theorem~\ref{BAMforRGB} that we are merely recomputing
the Gr\"obner basis~$\Sigma''$.
Moreover, we do not store the pairs of pairs, but the
$S$-vectors they generate, and we do not store $S$-vectors
which are clearly useless by part b) of the 
remark following Proposition~\ref{redSigma}.

The main difference occurs in steps~5) through~10).
Instead of computing the normal remainder of the S-vector,
we perform leading term reductions only and check the result
after each reduction step.
When we choose an element~$\sigma_{ij}$ in step~4),
it is not contained in~$\Cal{A}_d$, since if an element 
$\sigma_{k\ell}$ is appended to~$\Cal{A}$ in step~11)
and cannot be contained in~$\Cal{B}_d^\ast$ by step~10).
But the element~$\sigma_{ij}$
could have a leading term in~$\LT_\sigma(\Cal{A}_d)$
without being contained in~$\Cal{A}_d$. We claim that,
in this case, we know $\sigma_{ij} \TTo{\Cal{A}} 0$,
i.e.\ that~$\sigma_{ij}$ produces no new Gr\"obner
basis element.

To prove this claim, we first note that clearly~$\Cal{A}$
is a subtuple of~$\Cal{W}$ at all times.
Since the elements of~$\Cal{W}$ are fully interreduced,
the tail of~$\sigma_{ij}$ cannot be a leading term
of an element of~$\Cal{A}_d$. On the other hand, if
$\LT_\tau(\sigma_{ij})=\LT_\tau(\sigma_{i'j})$
for $\sigma_{i'j}\in \Cal{A}_d$, then the leading term
of the result of the reduction of~$\sigma_{ij}$
by~$\sigma_{i'j}$ is the tail of~$\sigma_{ij}$.
Hence~$\sigma_{ij}$ can be tail reduced using~$\Cal{A}_{<d}$.
By applying the same argument to the result of this
tail reduction step, we conclude that after several
tail reductions using~$\Cal{A}_{<d}$, we reach an element
of~$\Cal{A}_d$, and the claim follows.

The next possibility for~$\sigma_{ij}$
is that it is head irreducible with respect
to~$\Cal{A}$. In this case its leading term is
equal to $\LT_\tau(\sigma_{i'j})$ for some
$\sigma_{i'j}\in \Cal{W}_d$. Now Theorem~\ref{BAMforRGB}
says that we should put~$\NR_{\tau,\Cal{A}}(\sigma_{ij})$
into~$\Cal{A}$ and later swap it for~$\sigma_{i'j}$.
But, as we just saw, we can tail reduce~$\sigma_{ij}$
using~$\Cal{A}_{<d}$ until we reach~$\sigma_{i'j}$.
Thus the normal remainder is~$\sigma_{i'j}$
and is put into~$\Cal{A}$ immediately, i.e.\ without
actually performing the tail reductions and without a later
swap.

The last possibility for~$\LT_\tau(\sigma_{ij})$ is that
it can be reduced using~$\Cal{A}_{<d}$. This reduction step
is performed in step~7). Let us discuss the possible outcomes.

If the result is of the form $\tilde t\,\sigma_{k\ell}$
with $\tilde t\in\mathbb{T}^n\setminus \{1\}$, then~$\sigma_{k\ell}$
has a lower degree and satisfies $\sigma_{k\ell}\TTo{\Cal{A}}0$,
because~$\Cal{A}$ contains a truncated Gr\"obner basis.
Consequently, we have $\sigma_{ij}\TTo{\Cal{A}}0$
and step~4) of~\ref{BAMforRGB} tells us to try the next
S-vector.

If the result of the head reduction step has one of the
new leading terms provided by the elements of~$\Cal{W}_d$,
we notice this in step~8) or~9). In the first case,
the element of~$\Cal{V}_d$ has already been swapped
into~$\Cal{A}$ and nothing needs to be done. In the second
case, we perform the swap in step~9).

If the result is an element $\sigma_{k\ell}$ of degree~$d$
which can be further head reduced, we check in step~10)
whether $\sigma_{k\ell}\in \Cal{B}^\ast_d$. In that case
$\sigma_{ij}$ and~$\sigma_{k\ell}$ have the same
reductions and it suffices to treat~$\sigma_{k\ell}$
in step~7). Otherwise, we claim that~$\sigma_{k\ell}$
is one of the elements of~$B^\ast_d$ which has been
dealt with already, i.e.\ that we can go back to
step~4) and treat the next element of~$\Cal{B}_d^\ast$.

To prove this claim, we first use $\sigma_{ij}\in
\Cal{B}^\ast_d$ in order to write $\sigma_{ij}= t'\,
\sigma_{im} + t''\,\sigma_{jm}$ with $t',t''\in\mathbb{T}^n
\setminus \{1\}$ and $j<m\le s$.
Secondly, by step~7), we have the equality $\sigma_{ij}=
t\,\sigma_{i'j} \pm \sigma_{k\ell}$, where $\sigma_{k\ell}
=\pm \sigma_{i'i}$ and $t\in \mathbb{T}^n\setminus\{1\}$.
By looking at the coefficient of~$e_j$ in the
equation $\sigma_{ii'}=t'\, \sigma_{im} -t\,
\sigma_{i'j} - t''\, \sigma_{jm}$, we see that
$t\,\lcm(t_{i'},t_j)= t''\, \lcm(t_j,t_m)$. This term
is a multiple of~$t_{i'}$ and of~$t_m$. Hence it is of the
form $\tilde t\, \lcm(t_{i'},t_m)$ for some $\tilde t\in\mathbb{T}^n$,
and we have $\sigma_{ii'}=t'\,\sigma_{im} - \tilde t\,
\sigma_{i'm}$. If~$\tilde t\ne 1$, then
$\sigma_{k\ell}$ is a pair of pairs, i.e.\ it is either
in~$B_d^\ast$ or it is one of the elements of~$B_d^\ast$
treated before. Hence the claim follows if we can show that
$\tilde t=1$ does not happen.

Suppose that $\tilde t=1$. Then we are in the situation
of the lemma. Since the conditions of steps~8) and~9) did
not apply, it follows that~$\sigma_{k\ell}$ can be further
head reduced using~$\Cal{A}_{<d}$. Hence there exist
$u,u'\in\mathbb{T}^n$ and $j'<\max\{i,i'\}$ such that
$\sigma_{i'i}= u\,\sigma_{i'j'} + u'\,\sigma_{j'i}$
and $u\ne 1$ or $u'\ne 1$, depending on whether $i>i'$
or $i<i'$.

Now we show that $u'\ne 1$ is impossible. We use the
notation of the lemma and let $\delta=\deg_{x_\kappa}
(t_{j'})$, where~$x_\kappa$ is one of the indeterminates
occurring in~$t$, i.e.\ where case~1) of the lemma holds.
Then the equation $\lcm(t_{i'},t_i)= u\,\lcm(t_{i'},t_{j'})
=u'\, \lcm(t_i,t_{j'})$ shows $\max\{ \alpha,\alpha'\}
> \max\{\alpha,\delta\}$. This implies $\alpha' > \alpha$ and
$\alpha'>\delta$, in contradiction to case~1) of the lemma.
Similarly, we can show that $u\ne 1$ is impossible.
This concludes the proof of the claim.

Altogether, it follows that steps~5) -- 10) implement the full
reduction of~$\sigma_{ij}$ together with the swapping procedure of
step~7) of~\ref{BAMforRGB}. Hence the remaining elements
of~$\Cal{W}_d$ are precisely the minimal generators of degree~$d$
we are looking for, and they have to be appended to~$\Theta$ in
step~11).
\endproof

Let us apply this algorithm in the situation of
Example~\ref{ex1}.

\begin{example}\label{ex1cont}
Our task is to minimalize $\Cal{W}=\Sigma''=
\{\sigma_{12},\sigma_{13},\sigma_{24},\sigma_{34}\}$,
where we have $\deg_W(\sigma_{12})=9$,
$\deg_W(\sigma_{13})=10$, and $\deg_W(\sigma_{24})=
\deg_W(\sigma_{34})=9$.

In step~2), the algorithm constructs the set~$\Cal{B}^\ast$.
The pair of pairs $((2,4),(3,4))$ yields $S_{((2,4),(3,4))}
= z \sigma_{24} - x \sigma_{34} = -yz^2 \epsilon_2 + x^3 \epsilon_3 =
\sigma_{23}$, and this is the only element of~$\Cal{B}^\ast$.
Notice that it has degree~10.

In step~3), the algorithm starts to operate in degree $d=9$.
Since $\Cal{B}^\ast_9=\emptyset$, it appends $\sigma_{12}$,
$\sigma_{24}$, and $\sigma_{34}$ to~$\Cal{A}$
and~$\Theta$ in step~11).

Next we process degree~10. In step~4), we choose
$\sigma_{23}\in \Cal{B}^\ast_{10}$ and set $\Cal{B}^\ast_{10}=
\emptyset$. Then, in step~6), we find $\LT_\tau(
\sigma_{23})= x^3 \epsilon_3 = \LT_\tau (\sigma_{13})$,
where $\sigma_{13}\in \Cal{W}_{10}$. Hence~$\sigma_{13}$
is removed from~$\Cal{W}_{10}$ and appended to~$\Cal{A}$
in step~6).

Thus we have $\Cal{B}^\ast=\emptyset$ and~$\Cal{W}=
\emptyset$ at this point, and step~12) returns the set $\Theta=\{\sigma_{12},
\allowbreak \sigma_{24},\allowbreak \sigma_{34}\}$.
We note that this is the correct answer,
and there is an improvement over the application of
Proposition~\ref{GMrules} coming from the fact that in
step~6) we merely check
$\LT_\tau(\sigma_{ij})\in \LT_\tau(\Cal{W}_d)$ rather than
$\sigma_{ij}\in \Cal{W}_d$.
\end{example}

The following example provides a case where it is actually
necessary to do one head reduction step in~7) in order
to find a previously undiscovered non-minimal critical
syzygy.

\begin{example}\label{ex2}
Let $P=\mathbb{Q}[x_1,\dots,x_5]$ be standard graded, let
$r=1$ and $s=4$. The terms $t_1=x_2^2 x_3^6 x_4 x_5^2$,
$t_2= x_1^8 x_2 x_4 x_5^4$, $t_3= x_1^8 x_2^2 x_3^6$,
and $t_4=x_1^8 x_3^6 x_5^4$ yield the critical
syzygies $\sigma_{12}=x_1^8 x_5^2 \epsilon_1 - x_2 x_3^6 \epsilon_2$,
$\sigma_{13}= x_1^8 \epsilon_1 - x_4 x_5^2 \epsilon_3$,
$\sigma_{14}= x_1^8 x_5^2 \epsilon_1 - x_2^2 x_4 \epsilon_4$,
$\sigma_{23}= x_2 x_3^6 \epsilon_2 - x_4 x_5^4 \epsilon_3$,
$\sigma_{24}= x_3^6 \epsilon_2 - x_2 x_4 \epsilon_4$, and
$\sigma_{34}= x_5^4 \epsilon_3 - x_2^2 \epsilon_4$.
Here steps~1) and~2) of Proposition~\ref{GMrules}
discard~$\sigma_{23}$ and~$\sigma_{14}$,
because we have $\LT_\tau(\sigma_{23})=
x_4 x_5^4 \epsilon_3 = x_5^2 \LT_\tau(\sigma_{13})$
and $\LT_\tau(\sigma_{14})= x_2^2 x_4 \epsilon_4 =
x_2 \LT_\tau(\sigma_{24})$. Thus we have
$\Sigma''= \{ \sigma_{12},\mathstrut\sigma_{13},
\sigma_{24}, \sigma_{34}\}$. We note that we have $\deg_W(
\sigma_{12})=21$, $\deg_W(\sigma_{13})= 19$, and
$\deg_W(\sigma_{24}) = \deg_W(\sigma_{34}) = 20$.
But~$\mathstrut\Sigma''$ is not minimal, since
we have $\sigma_{12}=x_5^2 \sigma_{13} - x_2
\sigma_{24} +x_4 \sigma_{34}$.

Now we apply our algorithm. In step~2), we have to compute
$S_{((2,4),(3,4))} = x_2 \sigma_{24} - x_4 \sigma_{34}
= x_2 x_3^6 \epsilon_2 - x_4 x_5^4 \epsilon_3 = \sigma_{23}$.
Thus~$\sigma_{23}$ is appended
to~$\Cal{B}^\ast$. It has degree $\deg_W(\sigma_{23})=21$.
No further pairs of pairs are found.

In step~3), the algorithm starts to operate in
degree~$d=19$. We have $\Cal{B}^\ast_{19}=\emptyset$ and
$\Cal{W}_{19}=(\sigma_{13})$. Thus we append~$\sigma_{13}$
to~$\Cal{A}$ and $\Theta$ in step~11).
Next we pass to degree $d=20$.
We still have~$\Cal{B}^\ast_{20}=\emptyset$,
but now we get $\Cal{W}_{20}= (\sigma_{24},\sigma_{34})$.
In step~11), $\sigma_{24}$ and $\sigma_{34}$~are put
into~$\Cal{A}$ and~$\Theta$.

When we start processing degree~$d=21$, we have to choose
$\sigma_{23}\in \Cal{B}^\ast_{21}$ and set
$\Cal{B}^\ast_{21}=\emptyset$ in step~4). The leading term
$\LT_\tau(\sigma_{23})=x_4 x_5^4 \epsilon_3$ is not equal to
one of the leading terms of the elements of
$\Cal{A}_{21}$ or~$\Cal{W}_{21}$.
But we can perform a head reduction
step in~7), namely $\sigma_{23} \TTo{\sigma_{13}} -\sigma_{12}$. 
Here step~8) does not apply, but in step~9) we have
$\LT_\tau(\sigma_{12})\in\LT_\tau(\Cal{W}_{21})$.
Thus we continue by removing $\sigma_{12}$
from~$\Cal{W}_{21}$ and appending it to~$\Cal{A}$.

Finally, we get $\Cal{B}^\ast=\emptyset$ and $\Cal{W}=
\emptyset$. The algorithm returns
$\Theta= \{ \sigma_{13},\sigma_{24},\sigma_{34}\}$.
As mentioned above, the non-minimal critical
syzygy~$\sigma_{12}$ was discovered after one head reduction
step in~7).
\end{example}

\section{An Optimized Buchberger Algorithm}
\label{An Optimized Buchberger Algorithm}

In this section we combine the results obtained so far.
We continue to use the notation and conventions
of the previous sections. In particular, we let
\hbox{$P=K[x_1,\dots,x_n]$} be a polynomial ring over a
field~$K$ which is positively graded by a matrix~$W\in\Mat_{m,n}
(\mathbb{Z})$, and we let~$M$ be a graded submodule of a graded
free $P$-module $F=\dsum_{i=1}^r P(-\delta_i)$
which is generated by a tuple \hbox{$\Cal{V}=(v_1,\dots,v_s)$}
of homogeneous vectors. Furthermore, we 
let~$\sigma$ be a module term ordering on~$\Tnr$.

In the following theorem the sets of critical pairs
corresponding to the sets of critical syzygies considered
earlier are denoted by the normal letters corresponding 
to their calligraphic versions.

\begin{thm}{\sem (Optimized Buchberger Algorithm)}\\
\label{optBA}%
In the above situation, consider the following
sequence of instructions.

\begin{items}
\item[1)] Let $\Cal{W}=\Cal{V}$, $A=\emptyset$, $B=\emptyset$,
$B^\ast=\emptyset$, $\Cal{G}=\emptyset$, and let~$s'=0$.

\item[2)] Let~$d$ be the smallest degree
w.r.t.~{\tt Lex} of an element of~$B$
or~$\Cal{W}$. Form $B_d$, $B^\ast_d$, $\Cal{W}_d$, and delete their
entries from~$B$, $B^\ast$, and~$\Cal{W}$, respectively.

\item[3)] Apply ${\tt MinPairs}(A,B_d,B^\ast_d)$.

\item[4)] If $B_d=\emptyset$, then continue with step~7).
Otherwise, choose a pair $(i,j)$ in~$B_d$, delete it
from~$B_d$, and append it to~$A$.

\item[5)] Compute $S_{ij}$ and $S_{ij}'=
\NR_{\sigma,\Cal{G}}(S_{ij})$. If $S_{ij}'=0$, then
continue with~4).

\item[6)] Increase $s'$ by one, append $g_{s'}=S_{ij}'$
to~$\Cal{G}$, perform ${\tt Update}(B,B^\ast,g_{s'})$,
and continue with step~4).

\item[7)] If $\Cal{W}_d=\emptyset$ then continue with~10).
Otherwise, choose $v\in\Cal{W}_d$ and delete it in~$\Cal{W}_d$.

\item[8)] Compute $v'=\NR_{\sigma,\Cal{G}}(v)$. If $v'=0$,
continue with step~7).

\item[9)] Increase $s'$ by one, append $g_{s'}=v'$ to~$\Cal{G}$
and perform ${\tt Update}(B,B^\ast,g_{s'})$. Then continue with
step~7).

\item[10)] If $B=\emptyset$ and $\Cal{W}=\emptyset$,
then return $\Cal{G}$ and stop. Otherwise, continue with step~2).
\end{items}

\noindent Here the procedure ${\tt Update}(B,B^\ast,g_{s'})$
is defined as follows.

\begin{items}
\item[U1)] Form the set $C=\{(i,s') \mid 1\le i < s',\,
\gamma_i=\gamma_{s'}\}$.

\item[U2)] Delete from~$C$ all pairs $(j,s')$ for which there
exists an index $i$ in the set $\{1,\dots,j-1\}$ such that $t_{s'i}$
divides $t_{s'j}$.

\item[U3)] Delete from~$C$ all pairs $(i,s')$ for which
there exists an index $j$ in the set $\{i+1,\dots,s'-1\}$ such that
$t_{s'j}$ properly divides $t_{s'i}$.

\item[U4)] Find in~$C$ all pairs~$(i,s')$ and~$(j,s')$
such that $1\le i<j<s'$ and such that $\gcd(t_{is'},t_{js'})=1$.
For each of these, check if $(i,j)$ is already contained
in~$B^\ast$ and append it if necessary.

\item[U5)] Append the elements of~$C$ to~$B$ and stop.

\end{items}

Furthermore, the procedure ${\tt MinPairs}(A,B_d,B_d^\ast)$
is defined as follows.

\begin{items}
\item[M1)] If $B_d^\ast=\emptyset$, then
stop. Otherwise, choose a pair $(i,j)$
in~$B_d^\ast$ and remove it from~$B_d^\ast$.

\item[M2)] If $t_{ji}=t_{ji'}$ for some pair
$(i',j)\in A$, then continue with step~M1).

\item[M3)] If $t_{ji}=t_{ji'}$ for some pair
$(i',j)\in B_d$, then remove this pair
from~$B_d$ and append it to~$A$.
Continue with step~M1).

\item[M4)] Find $(i',j)\in A$ such that
$t_{ji'}$ divides~$t_{ji}$. Let $k=\min\{i,i'\}$,
and let $\ell=\max\{i,i'\}$. If $\gcd(t_{ij},t_{i'j})\ne
1$, then continue with~M1).

\item[M5)] If $t_{\ell k}=t_{\ell k'}$ for some
pair $(k',\ell)\in A$, then continue with M1).

\item[M6)] If $t_{\ell k}=t_{\ell k'}$ for some
pair $(k',\ell)\in B_d$, then delete this pair in~$B_d$,
append it to~$A$, and continue with~M1).

\item[M7)] If $(k,\ell)\in B_d^\ast$, then delete $(k,\ell)$
in~$B_d^\ast$ and continue with~M4), applied to this pair.

\item[M8)] Continue with step~M1).

\end{items}

Altogether, we obtain an algorithm which computes a tuple
$\Cal{G}$ whose elements form a homogeneous
$\sigma$-Gr\"obner basis of~$M$. 
Moreover, the set of pairs which
are treated at some time in steps 4) -- 6) of the algorithm
corresponds to a minimal system of generators of the module
$\Syz_P(c_1t_1e_{\gamma_1},\dots,c_{s'}t_{s'}e_{\gamma_{s'}})$.

\end{thm}

\proof
The main algorithm of this theorem agrees with the
Homogeneous Buchberger Algorithm (see Theorem~\ref{homBA}), except
for the introduction of the procedure ${\tt MinPairs}
(A,B_d,B_d^\ast)$ in step~3) and the alteration of
the enlargement of~$B$ in steps~5) and~8) of Theorem~\ref{homBA}
which is now performed by the procedure ${\tt Update}
(B,B^\ast,g_{s'})$.

The foundation for these changes is the material presented above,
especially Proposition~\ref{SigmaBAM}. In steps 4) -- 6) we want 
to treat only those pairs $(i,j)$ for which the corresponding
elements~$\sigma_{ij}$ are contained in the minimal system of
generators~$\Theta$ of the graded $P$-module
$\Syz_P(c_1t_1e_{\gamma_1},\dots,c_{s'}t_{s'}e_{\gamma_{s'}})$.

Procedure ${\tt Update}(B,B^\ast,g_{s'})$ applies Rules~1) and~2)
of Gebauer-M\"oller in steps~U2) and~U3), respectively. Moreover,
notice that step~U4) computes all pairs of pairs which satisfy the
condition of  part b) of Remark~\ref{remBAMapl}, and stores the pairs
corresponding to the resulting S-vectors in~$B^\ast$.

Thus, in order to minimalize the critical pairs we process,
we need to apply Proposition~\ref{SigmaBAM} to the
set of critical syzygies corresponding to the
set of critical pairs~$B$, where we can refrain from
computing the pairs of pairs, because they have already
been generated and stored in~$B^\ast$. This task
is performed by the procedure ${\tt MinPairs}(A,B_d,B_d^\ast)$.
Its steps~M1) --~M8) are easy translations of steps~4) --~10)
of Proposition~\ref{SigmaBAM} into the language of pairs.
Notice that we have $\LT_\tau(\sigma_{ij})=\LT_\tau
(\sigma_{k\ell})$ if and only if $j=\ell$ and $t_{ji}=
t_{\ell k}$.
Altogether, ${\tt Update}(B,B^\ast,g_{s'})$ and ${\tt
MinPairs}(A,B_d,B_d^\ast)$ make sure that only the pairs
corresponding to~$\Theta$ are treated at some point in steps~3)
--~6). 

Finally, we remark that~$A$ is used to keep track of the
pairs $(i,j)$ for which $\sigma_{ij}$ is in that part of the
minimal $\tau$-Gr\"obner basis~$\Sigma''$
of~$\Syz_P(c_1t_1e_{\gamma_1},\dots,c_{s'}t_{s'}e_{\gamma_{s'}})$
which has been computed so far. Thus it is updated when a
non-minimal element of~$\Sigma''$ is found in step~M3) or step~M6),
and when a pair corresponding to an element of~$\Theta$ is chosen
for treatment in step~3).
\endproof

Let us illustrate the performance of this algorithm
by a simple example. It shows that cases like Example~\ref{ex1}
occur naturally during actual Gr\"obner basis computations.

\begin{example}\label{ex3}
Let $P=\mathbb{Q}[x,y,z]$ be standard graded, let 
$\sigma={\tt DegLex}$, let \hbox{$r=1$}, and let
$M\subseteq P$ be the homogeneous ideal generated by
the polynomials $v_1=x^3z^2+x^2y^2z$, $v_2=x^3y^8$, and $v_3=y^{10}z^2$.
Then the leading terms are $t_1=x^3z^2$, $t_2=x^3y^8$,
and $t_3=y^{10}z^2$. Let us follow the steps of the
Optimized Buchberger Algorithm.

The first degree is~$d=5$. Since $B_5=\emptyset$,
the first actions are to choose $v_1\in\Cal{W}_5$
in step~7) and append $g_1=v_1$ to~$\Cal{G}$ in step~9).
Then we continue with $d=11$ and choose
$v_2\in\Cal{W}_{11}$ in step~7). Since $v'=\NR_{
\sigma,\Cal{G}}(v_2)=v_2$, we append$g_2=v_2$ to~$\Cal{G}$
in step~9) and update the set of pairs. The
result is $B=\{(1,2)\}$ and $B^\ast=\emptyset$.
Now we have to treat the degree $d=12$.
Notice that the degree of the pair $(1,2)$ is~13.
Hence $B_{12}=\emptyset$ and we have to choose
$v_3\in\Cal{W}_{12}$ in step~7). Since $v'=\NR_{
\sigma,\Cal{G}}(v_3)=v_3$, we append$g_3=v_3$ to~$\Cal{G}$
in step~9) and update the set of pairs.
In step~U1), we form $C=\{(1,3),(2,3)\}$. In step~U2),
we obtain $t_{31}=x^3=t_{32}$, and therefore $(2,3)$
is deleted in~$C$. The result is $B=\{(1,2),(1,3)\}$
and $B^\ast=\emptyset$. This completes degree~12, and
we continue with degree~13.

We choose the pair $(1,2)$ in step~4) and append it
to~$A$. Then we compute $S_{12}=y^8 g_1 -z^2 g_2 =
x^2y^{10}z$ and $S'_{12}=\NR_{\sigma,\Cal{G}}(S_{12})
=x^2 y^{10} z$. Thus we have a new Gr\"obner basis
element $g_4= x^2 y^{10} z$ and need to update
the pairs again. In step~U1), we form $C=\{(1,4),
(2,4), (3,4)\}$. Step~U2) does not apply, but in step~U3)
we remove the pair $(1,4)$ from~$C$, since
$t_{42}=x$ properly divides $t_{41}=xz$.
Now we check that $t_{24}=y^2z$ properly divides
$t_{23}=y^2z^2$ and $t_{34}=x^2$ properly divides
$t_{32}=x^3$. Hence the pair $(2,3)$ is appended
to~$B^\ast$.

At this point we have finished degree~13, and we have
the following situation: $A=\{(1,2)\}$, $B=\{(1,3), (2,4),
(3,4)\}$, $B^\ast = \{(2,3)\}$,
$\Cal{G}=\{g_1,\dots,g_4\}$, and $s'=4$. The next degree is~$d=14$,
where we have to deal with the pairs in $B_{14}=
\{(2,4),(3,4)\}$. Since $B^\ast_{14}=\emptyset$,
we choose $(2,4)$ in step~4) and append it to~$A$.
Then we compute $S_{24}=0$ and continue by choosing
$(3,4)$ in~$B_{14}$ and adding it to~$A$.
Again $S_{34}=0$, and degree~14 is finished.

Now we start degree~15 by performing ${\tt MinPairs}
(A,B_{15},B_{15}^\ast)$, where we have $A=\{(1,2),\allowbreak
(2,4),\allowbreak (3,4)\}$,
$B_{15}=\{(1,3)\}$ and $B^\ast_{15}=\{(2,3)\}$.
In step~M1), we choose $(2,3)$. In step~M3),
we discover $t_{32}=x^3=t_{31}$, where $(3,1)\in B_{15}$.
Hence $(1,3)$ is deleted in~$B_{15}$ and appended to~$A$.
Then the procedure is finished, and the facts that
$B_{15}=\emptyset$ as well as $\Cal{W}_{15}=\emptyset$
allow us to return $\Cal{G}$ and stop.

As in Example~\ref{ex1}, we have found one useless pair, namely the
pair $(1,3)$ in degree~15, which would not have been discovered by
the Gebauer-M\"oller Installation, and which we were able to
discard by a simple combinatorial check.
\end{example}

\begin{remark} Let us discuss the efficiency of
the algorithm of Theorem~\ref{optBA}.

\begin{items}
\item Steps~U2) and~U3) of this algorithm correspond to
Rules~1) and~2) of the Gebauer-M\"oller installation. However,
Rule~3 is not performed by the procedure ${\tt Update}(\dots)$, but
by step~M2) of the procedure ${\tt MinPairs}(\dots)$. In fact,
step~M2) gets rid of more pairs than Rule~3, because Rule~3
requires $(i,j)\in B_d^\ast
\cap B_d$, whereas we only need a pair $(i,j)
\in B_d^\ast$ such that $\LT_\tau(\sigma_{ij})=\LT_\tau(\sigma_{i'j})$
for some $(i',j)\in B_d$.

\item A potential drawback of our approach is that
the number of pairs of pairs considered in step~U4)
is quadratic in the number of elements of~$C$ surviving
steps~U2) and U3). But that number is usually fairly small.
Hence the cost of~U4) and the cardinality of~$B^\ast$
tend to be rather small.
On the other hand, we do not need to check Rule~3 for all
elements of the list~$B$ which is usually rather long.
Our experiments suggest that, on average, the overhead of the
two approaches is comparable.

\item Our procedure ${\tt MinPairs}(\dots)$ is very efficient
in treating the elements of~$B_d^\ast$. Each time we loop through
steps M2) -- M8), we delete one pair in~$B_d^\ast$, and~$B_d^\ast$
is never enlarged. In practice, we find that the lists~$B_d^\ast$
are generally small. Hence our algorithm harnesses the full power
and efficiency of the Gebauer-M\"oller installation, while it
simultaneously kills {\em all} unnecessary pairs at a comparatively
small cost.

\end{items}
\end{remark}

\section{Experimental Data and Conclusions}
\label{Experimental Data and Conclusions}

In this section we want to provide the reader with some
experimental numerical data which illustrate the performance of the
Optimized Buchberger Algorithm~\ref{optBA} as well as technical
observations coming from an implementation in an experimental
version of the ``\cocoa\ 5'' library in C++.

In the following table, we compare the application of Rules 1)
--~3) of~Proposition~\ref{GMrules} to our procedures ${\tt Update}(\dots)$
and ${\tt Minpairs}(\dots)$ in Theorem~\ref{optBA}, i.e.\ to the
algorithm of Proposition~\ref{SigmaBAM}. 
Let us point out that our procedure always minimalizes the critical
pairs, independent of the order of the underlying terms.
(Non-minimal critical pairs are recognized at different steps,
though.) For the Gebauer-M\"oller installation, however, the number
of undiscovered non-minimal critical pairs depends strongly on this
order.

To aid the reader in understanding this table, let us explain the meaning
of the symbols.

\begin{items}
\item[$\bullet$] $\#(G)$ is the cardinality of the reduced Gr\"obner basis
of the corresponding ideal.

\item[$\bullet$] $\#(\Sigma)$ is the total number of pairs, i.e.\
$\#(\Sigma) = {\#(G)\choose 2}$.

\item[$\bullet$] $\#(\Sigma'')$ is the number of pairs surviving Rules~1)
and~2), i.e.\ the cardinality of the reduced Gr\"obner basis of pairs.

\item[$\bullet$] $B$ is the number of pairs killed by Rule 3),
the Gebauer-M\"oller ``Backwards'' criterion.

\item[$\bullet$] $M23$ is the number of pairs killed by steps~M2) and~M3) in
Theorem~\ref{optBA}.

\item[$\bullet$] $M48$ is the number of pairs killed by steps~M4) --~M8)
in Theorem~\ref{optBA}.

\item[$\bullet$] ${\bf Gain}=M23+ M48-B$,
i.e.\ the number of newly discovered
non-minimal critical pairs.

\item[$\bullet$] $\#(\Theta)$ is the cardinality of a minimal system of
generators of the syzygies of the leading terms. Hence we have
$\#(\Theta) = \#(\Sigma'') - M23 - M48$
\end{items}

\bigskip\bigbreak

\def\newrow{\cr\noalign{\hrule}}
\def\Newrow{\cr\noalign{\hrule height 1pt}}
\def\vr{\vrule}
\def\qq{\hskip 0.05em}
\def\q{\hskip 0.05em}

{\small

\centerline{\vbox{\halign{\strut
  \vrule width 1pt height 11pt\q # \hfil
&\vr\ \vr  \qq \hfil  # \q
&\vr       \qq \hfil  # \q
&\vr       \qq \hfil  # \q
&\vr       \qq \hfil  # \q
&\vr       \qq \hfil  # \q
&\vr       \qq \hfil  # \q
&\vr       \qq \hfil  # \q
&\vr\   \qq \hfil  # \q  \vrule width 1pt height 11 pt\cr
\noalign{\hrule height 1pt}
 & $\#(G)$ \  & $\#(\Sigma)\ $\ &$ \#(\Sigma'') $\ & $ B $\  &$ M23 $\
 & $ M48 $\ & {\bf Gain}\ & $\#(\Theta)$\
 \newrow
T\^\relax 51      & 83 & 3,403    & 250   &  7 &  7   & 0 &${\mathbf{0}}$  & 243
\newrow
Twomat3      & 109 & 5,886   & 741   &  15 &  26   & 1 &${\mathbf{12}}$  & 714
\newrow
Alex3      & 211 & 22,155  & 684   & 54   & 56    & 1 & ${\mathbf{3}}$ &  627
\newrow
Gaukwa4    & 267 & 35,511  & 1,772 & 101  & 113  & 3 & ${\mathbf{15}}$ & 1,656
\newrow
Kin1       & 306 & 46,665  & 3,411 & 70  & 172  & 0 & ${\mathbf{102}}$ & 3,239
\newrow
Wang (Lex) & 317 & 50,086  & 1,457 & 60  & 61   & 7 & ${\mathbf{8}}$ & 1,389
\newrow
Cyclic 7   & 443 & 97,903   & 2,651   & 17     & 17     & 0 & ${\mathbf{0}}$ 
& 681
\newrow
Hairer-2   & 506 & 127,765 & 5,305 & 150  & 152  & 4 & ${\mathbf{6}}$ & 5,149
\newrow
Hom-Gonnet & 854 & 364,231 & 11,763  & 587 &  648  & 27 & ${\mathbf{88}}$ 
&11,088
\newrow
Mora-9     & 4,131 & 8,530,515 & 46,395  &  1,930  & 1,914  & 23 & 
${\mathbf{7}}$ & 44,458
\Newrow
}}}

}

\bigskip

The rows of this table correspond to standard examples of Gr\"obner
basis computations. A file containing a description of every example
can be downloaded~at

\begin{center}
{\tt ftp://cocoa.dima.unige.it/papers/CaboaraKreuzerRobbiano03.cocoa}
\end{center}

Moreover, a file containing the list of leading terms of the reduced 
Gr\"obner basis for each example can be downloaded at
\medskip

\hbox{{\tt ftp:/\hskip -2pt 
/cocoa.dima.unige.it/papers/CaboaraKreuzerRobbiano03\_2.cocoa}}

\noindent {\it Technical note:}\/ In the well-known example ``Cyclic 7" 
we have homogenized using a new {\it smallest}\/ 
indeterminate (see the file mentioned above).

\smallskip

For the reader who would like to run his own tests, we note that
$\#(G)$, $\#(\Sigma)$, and $\#(\Theta)$ are invariants of the reduced
Gr\"obner basis. But the effect of both the Gebauer-M\"oller
installation and our Optimized Buchberger Algorithm depend strongly on the
order in which the elements of~$\Sigma$ are produced during a
Gr\"obner basis computation. For instance, this means that it
depends on the chosen selection strategy. In our implementation pairs are
kept ordered in increasing {\tt DegLex} ordering, reductors are kept in 
the same order they are produced, reductors of the same degree are kept
interreduced, and the reduction strategy is full reduction.

The following table shows some timings. It compares Singular 2.0.0
with the current experimental version of \cocoa\ 5 using the GM 
and CKR pair handling algorithms. Timings are in seconds for Linux 
running on an Athlon 2000+ CPU with 1.5GB RAM. 
All computations are over the rationals where the timings of the
base field operations in Singular and \cocoa\ seem to be 
comparable.

\smallskip

\noindent {\it Technical note:}\/ The reason why we include a 
comparison with Singular is an explicit request made
by a referee, who suggested comparing our timings with 
``another efficient implementation". The table below indicates
that both Singular and \cocoa\ 5 have efficient implementations
of the Buchberger algorithm, and that our new algorithm has at
least the same efficiency.

\bigskip

{\small
\centerline{\vbox{\halign{\strut
  \vrule width 1pt height 11pt\q # \hfil
&\vr\ \vr  \qq \hfil  # \q
&\vr       \qq \hfil  # \q
&\vr       \qq \hfil  # \q
\vrule width 1pt height 11 pt\cr
\noalign{\hrule height 1pt}
\qquad & Singular 2.0.0 \  & \cocoa 5 GM & \cocoa 5 CKR \newrow
T\^\relax 51 (Lex)   & 149.32 & 7.28    & 7.14 \newrow
Twomat3      & 1.21 & 8.66   & 8.50 \newrow
Alex3        & $<\!<$1 & 0.54  & 0.56 \newrow
Gaukwa4      & 80.30 & 99.31  & 98.57 \newrow
Kin1         & 407.09 & 89.25  & 87.41 \newrow
Wang (Lex)   & $>$1200 & 382.86  & 379.31 \newrow
Cyclic 7     & $>$1200 & 76.61   & 76.65 \newrow
Hairer-2     & 79.36 & 141.83 & 139.76 \newrow
Hom-Gonnet   & 3.97 & 4.55 & 4.95  \newrow
Mora-9       & 30.53 & 86.17 & 89.75  \Newrow
}}}

}

\bigskip\bigbreak

\begin{center}
{\bf Conclusions}
\end{center}

\medskip

First of all, let us collect some
technical observations based on our implementation of the Optimized
Buchberger Algorithm.

\begin{items}
\item When we apply Rules~1 and~2 of Proposition~\ref{GMrules},
the remaining set of pairs~$\Sigma''$ is usually almost a minimal
system of generators
of the module $\Syz_P(c_1t_1e_{\gamma_1},\dots,c_{s'}t_{s'}e_{\gamma_{s'}})$.
Thus both Rule~3 and our algorithm kill comparatively few pairs.
Nonetheless, over the rationals (or other costly fields), the
saving is worthwhile because the treatment of each single pair can
take a long time.

\item Steps~M5) -- M7) in the Optimized Buchberger Algorithm are
independent. Hence it is possible to order them in such a way
that the computational cost is minimized. This may be important
if there is a large number of elements in~$B_d^\ast$ to be processed,
since the operations may have substantially different computational
costs.

\item All operations in our procedures {\tt Update(...)} and
{\tt MinPairs(...)} have been greatly eased by memorizing the terms
$t_{ij}, t_{ji}$ and $\lcm(t_i,t_j)$ directly in the pair data type.

\item When a search is performed on the pairs in $A$, $B$, or $B_d$,
full advantage can be taken of the fact that we may
rely on data structures which allow logarithmic search costs.

\end{items}

Looking at the timings above, we see that, on average and with
comparable implementations, our new algorithm is faster 
than the Gebauer-M\"oller installation. In some examples, the gains
are relatively small, and in exceptional cases, the structure of the 
combinatorial data produces a larger overhead for our algorithm
than for the Gebauer-M\"oller installation.

\subsubsection*{Acknowledgments}

The third author is grateful to the organizers of the First
International Con\-gress of Mathematical Software (Beijing 2002)
for their hospitality and the possibility to present his
work to a wide audience. An extended abstract of this
paper was published in the proceedings of the congress
(see \citep{CKR2002}).

\newpage


\end{document}